\documentclass[11pt,onecoulme]{article}
\usepackage{amsfonts}
\usepackage{mathrsfs}
\usepackage{amssymb}
\usepackage{color, amsmath,amssymb, amsfonts, amstext,amsthm, latexsym}

\usepackage{amssymb, epsfig, amssymb, latexsym}
\usepackage{amsmath}
\usepackage{graphicx}
\usepackage{longtable}

\numberwithin{equation}{section}

\allowdisplaybreaks

\newtheorem{definition}{Definition}[section]
\newtheorem{theorem}{Theorem}[section]
\newtheorem{lemma}{Lemma}[section]
\newtheorem{remark}{Remark}[section]

\newtheorem{hyp}[definition]{Hypothesis}
\title{Viscosity Solutions to First Order Path-Dependent  Hamilton-Jacobi-Bellman Equations in Hilbert Space \thanks{This work was partially supported by  the National Natural Science Foundation of China  (Grant No. 11401474), Shaanxi Natural Science Foundation
               (Grant No. 2017JM1016) and the Fundamental Research Funds for the Central Universities (Grant No. 2452019075).}}

\author{  Jianjun Zhou   \\
       College of Science,
             Northwest A\&F University,\\ Yangling 712100, Shaanxi, P. R.
             China\\
      \emph{E-mail:zhoujianjun@nwsuaf.edu.cn} }
          \date{}
\begin{document}

\maketitle

\pagestyle{plain}

\begin{abstract}
In this article, a notion of viscosity solutions  is introduced for first order path-dependent  Hamilton-Jacobi-Bellman (PHJB)
          equations  associated with   optimal control problems for path-dependent   evolution equations in Hilbert space.
     We
             identify the value functional of  optimal  control problems as
             unique viscosity solution to the associated PHJB equations. We also show that our notion of viscosity
solutions is consistent with the corresponding notion of classical solutions,
and satisfies a stability property.
\medskip

\end{abstract}


\section{Introduction}
%%%%%%%%%%%%%%%%%%%%%%%%%%%%%%%%%%%%%%%%%%%%%%%%%%%%%%%%%%%%%%%%%

%The notion of  viscosity solutions for  Hamilton-Jacobi-Bellman (HJB) equations, first introduced in 1983 by Crandall and Lions \cite{cra1}, has  become an indispensable tool in optimal control theory and numerous subjects
%related to it. We refer to the survey paper of Crandall, Ishii and  Lions \cite{cran2} and  the monographs of Fleming
%      and Soner \cite{fle1} and Yong and Zhou \cite{yong} for a detailed account
%for the theory of viscosity solutions. For viscosity solutions in infinite dimensional Hilbert spaces, we refer to Gozzi,  Rouy and \'{S}wi\c{e}ch \cite{gozz3}, Lions \cite{lio1}, \cite{lio2}, \cite{lio3} and \'{S}wi\c{e}ch \cite{swi}.
%\par
Viscosity solutions for first order Hamilton-Jacobi-Bellman (HJB) equations in infinite dimensions  have been investigated by Crandall and Lions in \cite{cra3,cra4,cra5} for the case without unbounded term, in \cite{cra6,cra7,cra9} for the case with unbounded linear term,
and in \cite{cra8} for the case with unbounded nonlinear term. Soon after,  Gozzi, Sritharan, and \'{S}wi\c{e}ch \cite{gozz} studied Bellman equations associated to control problems for variational solutions of
the Navier-Stokes equation. We also mention the work of  Li and Yong \cite{li}, where the  general unbounded first order HJB equations in infinite dimensional Hilbert spaces are studied.
%\par
% For the path-dependent  case, the theory of viscosity solutions is more difficult.
% Luyakonov \cite{luk} developed
%a theory of viscosity solutions to fully non-linear path-dependent first order  Hamilton-Jacobi equations.
%The existence and uniqueness theorems
%are proved when Hamilton function $\mathbf{H}$ is $d_p$-locally Lipschitz continuous in the path function.
%For the stochastic path-dependent  case, the notion of viscosity solutions  was introduced by Ekren, Keller, Touzi
% and Zhang \cite{ekren1} in the semilinear context, and further extended to the fully nonlinear case by Ekren,
%Touzi and Zhang  \cite{ekren2}, \cite{ekren3} and \cite{ekren4}, Ekren \cite{ekren0} and Ren \cite{ren}.
%The uniqueness results for the fully nonlinear case are only valid if  the  Hamilton function $\mathbf{H}$ is uniformly nondegenerate.
%    Ren,  Touzi and  Zhang \cite{ren1}  studied  the degenerate case and established  the comparison principle  when  the Hamilton function  $\mathbf{H}$ is $d_p$-uniformly continuous with respect to  the path function.
\par
Fully nonlinear  path-dependent first order  Hamilton-Jacobi equations have been studied by Lukoyanov  \cite{luk}.  The existence and uniqueness theorems
are obtained  when Hamilton function $\mathbf{H}$ is $d_p$-locally Lipschitz continuous in the path function. In our paper \cite{zhou}, we extended the results in \cite{luk} to $d_\infty$-Lipschitz continuous case.
Bayraktar and Keller \cite{bay}  proposed notions of minimax and viscosity solutions for a
class of fully nonlinear path-dependent HJB (PHJB) equations with nonlinear,
monotone, and coercive operators on Hilbert space and proved
the existence, uniqueness and stability
for minimax solutions.
For the second order path-dependent case, a viscosity solution
approach has been successfully initiated by Ekren, Keller, Touzi
 and Zhang \cite{ekren1} in the semilinear context, and further extended to
 fully nonlinear
equations by Ekren, Touzi, and Zhang \cite{ekren3, ekren4},  elliptic
equations by Ren \cite{ren}, obstacle problems by Ekren \cite{ekren0}, and degenerate second-order
equations by Ren, Touzi, and Zhang \cite{ren1}.
Cosso, Federico, Gozzi, Rosestolato, and Touzi \cite{cosso} studied a class of semilinear second order PHJB equations with a linear unbounded operator on Hilbert
space.
\par
In this paper, we consider the following   controlled path-dependent evolution
                 equation (PEE):
\begin{eqnarray}\label{state1}
\begin{cases}
            \dot{X}^{\gamma_t,u}(s)=AX^{\gamma_t,u}(s)+
           F(X_s^{\gamma_t,u},u(s)),  \ \ s\in [t,T],\\
~~~~X_t^{\gamma_t,u}=\gamma_t\in {\Lambda}_t.
\end{cases}
\end{eqnarray}
In the above equation,  $\Lambda_t$ denotes the set of all
 continuous $H$-valued functions defined over $[0,t]$ and ${\Lambda}=\bigcup_{t\in[0,T]}{\Lambda}_{t}$;
                the unknown $X^{\gamma_t,u}(s)$, representing the state of the system, is an $H$-valued process;
                the control process $u$ takes values in some  metric space $(U,d)$ ; $A$ is the generator of
         a $C_0$  semigroup of bounded linear operator in Hilbert space
        $H$ and the coefficient $F$ is assumed to satisfy $d_\infty$-Lipschitz condition.% with
       % respect to appropriate norm. %Under suitable assumptions, there exists a unique
%                adapted process $X^{\gamma_t,u}(s),\ s\in[t,T]$, solution to
%                (\ref{state1}).
\par
              We try  to minimize a cost functional of the form:
\begin{eqnarray}\label{cost1}
                     J(\gamma_t,u):=\int_{t}^{T}q(X^{\gamma_t,u}_\sigma,u(\sigma))d\sigma
                          +\phi(X^{\gamma_t,u}_T),\ \ \ (t,\gamma_t)\in [0,T]\times {\Lambda},
\end{eqnarray}
 over all admissible controls ${{\mathcal
                  {U}}}[t,T]$.
Here  $q$ and $\phi$ are given real functionals on ${\Lambda}\times U$ and ${\Lambda}_T$, respectively.
             We define the value functional of the  optimal
                  control problem as follows:
\begin{eqnarray}\label{value1}
V(\gamma_t):=\inf_{u\in{\mathcal
                  {U}}[t,T]}J({\gamma_t,u}),\ \ \ (t,\gamma_t)\in [0,T]\times {\Lambda}.
\end{eqnarray}
  The goal of this article is to characterize this value functional  $V$. We  consider the following
                     PHJB equation:
  \begin{eqnarray}\label{hjb1}
\begin{cases}
\partial_tV(\gamma_t)+(A^*\partial_xV(\gamma_t),\gamma_t(t))_H+{\mathbf{H}}(\gamma_t,\partial_xV(\gamma_t))= 0,\ \ \  (t,\gamma_t)\in
                               [0,T)\times {\Lambda},\\
 V(\gamma_T)=\phi(\gamma_T), \ \ \ \gamma_T\in {\Lambda}_T;
 \end{cases}
\end{eqnarray}
      where
$$
                                {\mathbf{H}}(\gamma_t,p)=\sup_{u\in{
                                         {U}}}[
                        (p,F(\gamma_t,u))_H+q(\gamma_t,u)], \  \ (\gamma_t,p)\in {\Lambda}\times H.
$$
 $A^*$ denotes  the adjoint operator of $A$. The definitions
of
                         $ \partial_t$ and $\partial_{x}$ will be introduced in Section 2.
\par
                         In this paper we will
                         develop a concept  of  viscosity solutions
                         to
                         PHJB equations on the space of  $H$-valued continuous paths (see Definition \ref{definition4.1} for details) and show that the value functional
                         $V$  defined in  (\ref{value1}) is   unique viscosity solution to the PHJB equation (\ref{hjb1}).
\par
The main challenge %for our $H$-valued path-dependent case
comes from the both facts that the path space $\Lambda_T$ is an infinite dimensional Banach space in which the maximal norm $||\cdot||_0$ is not smooth, and the operator $A$ is unbounded. In \cite{zhou} we gave a definition of viscosity solutions in a sequence of
bounded and uniformly Lipschitz continuous functions spaces
$C^\mu_{M_0}$
which are compact subsets in $R^d$-valued path-dependent case, and proved that the value functional is unique viscosity solution to the associated PHJB equation. However, %the spaces %bounded and uniformly Lipschitz continuous functions spaces
$C^\mu_{M_0}$  are not compact in $H$-valued path-dependent case.
By studying the $B$-continuity of the value functional, Li and Yong \cite{li} proved  the value functional is unique viscosity solution to the associated HJB equation in Hilbert space under the assumption $(A.4)$ on page 231 of {\cite{li}}. However, the value  functional  does not have the $B$-continuity  under the assumption $(A.4)$ in path-dependent  case. Therefore, the  techniques introduced in \cite{zhou} and \cite{li} are not applicable in
                 our case.
\par
 The main contribution of this paper
is the introduction of
%two
an appropriate  functional $\overline{\Upsilon}^2(\cdot,\cdot)$ %and $\Upsilon^2(\cdot)$
on $\Lambda\times \Lambda$.
% and $\Lambda$, respectively.
%For every fixed $(t,\eta_t) \in [0,T)\times\Lambda$,
% $\Upsilon^2(\cdot,\eta_t)\in C^{1}(\Lambda^{t})$ (see Lemma \ref{theoremS}).
The  functional is the key to prove the
 stability and uniqueness of viscosity solutions. Using Lemma \ref{theoremS000},
  we can define an auxiliary function $\Psi$ which includes the  functional $\overline{\Upsilon}^2(\cdot,\cdot)$ (see Step 1 in the proof of Theorem \ref{theoremhjbm}).
  More importantly, we can use $\overline{\Upsilon}^2(\cdot,\cdot)$ to  define a  smooth  gauge-type function and apply
                 a modification of Borwein-Preiss variational principle (see Lemma \ref{theoremleft}) to get a maximum of a perturbation of the auxiliary function $\Psi$.
 %Thanks to Lemma \ref{theoremS000},
 Then we prove the uniqueness of viscosity solutions without  the assumption $(A.4)$ of \cite{li}.  Regarding existence, we  prove that the value functional  $V$  defined in  (\ref{value1}) is  a viscosity solution to  PHJB (\ref{hjb1})   under our definition by functional It$\hat{o}$  formula and dynamic programming principle.
% We
%also observe that, in contrast with the  case without path-dependent, our
%well-posedness result applies to equations which cannot be treated, up to now, with the known theory of
%viscosity solutions.
  We also emphasize that, with respect to the standard viscosity solution theory in infinite dimension, the  assumption $(A.4)$ of  \cite{li} is completely bypassed in our framework. Therefore, even in  the  case without path-dependent, our
well-posedness result applies to equations which cannot be treated, up to now, with the known theory of
viscosity solutions.

\par
             The paper is organized as follows. In the following
              section, we introduce  preliminary results on  path-dependent optimal control problems in Hilbert space,  and prove
          Theorem \ref{theoremito} and  Lemmas \ref{theoremS} and  \ref{theoremS000} which are the key of the existence and uniqueness results of viscosity solutions.
             In Section 3, we introduce our notion of viscosity solutions to
           equation (\ref{hjb1}) and  prove  that the value functional $V$ defined by (\ref{value1}) is a viscosity solution. We also show
             the consistency with the notion of classical solutions and the stability result.
                  In Section 4,
               the uniqueness of viscosity solutions for equation (\ref{hjb1}) is proven and  Section 5 is devoted to proving $(\hat{\Lambda}^t, d_{\infty})$ and $({\Lambda}^t, d_{\infty})$ are two complete metric spaces.

\section{Preliminary work}  \label{RDS}
%%%%%%%%%%%%%%%%%%%%%%%%%%%%%%%%%%%%%%%%%%%%%%%%%%%%%%%%%%%%%%%%%

\par
 We list some notations that are used in this paper. We use the symbol $|\cdot|$ to denote the norm in
         a Banach space $\Xi$, with a subscript if necessary. Let  $H$
         denote a real separable Hilbert space, with scalar products
         $(\cdot,\cdot)_H$. %Let
%         $L(H)$ denotes the space of all
%         bounded linear operators from $H$ into itself.
            The operator $A$ is the generator of a strongly continuous
                 semigroup $\{e^{tA}, t\geq0\}$ of bounded linear operators in the
                 Hilbert space $H$.
                 The domain of the
         operator $A$ is denoted by ${\mathcal {D}}(A)$. $A^*$ denotes  the adjoint operator of $A$.
 %The domain of a linear (unbounded)
%         operator $A$ is denoted ${\mathcal {D}}(A)$.
         Let % $n$ be a positive integer and
$T>0$ be a fixed number.  For each  $t\in[0,T]$,
          define
         $\hat{\Lambda}_t:=D([0,t];H)$ as  the set of c$\grave{a}$dl$\grave{a}$g  $H$-valued
         functions on $[0,t]$.
       We denote $\hat{\Lambda}^t=\bigcup_{s\in[t,T]}\hat{\Lambda}_{s}$  and let  $\hat{\Lambda}$ denote $\hat{\Lambda}^0$.
       \par
%A very important remark on the notation $\hat{\Lambda}$:
As in  \cite{dupire1}, we will denote elements of $\hat{\Lambda}$ by lower
case letters and often the final time of its domain will be subscripted, e.g. $\gamma\in \hat{\Lambda}_t\subset \hat{\Lambda}$ will be
denoted by $\gamma_t$. Note that, for any $\gamma\in  \hat{\Lambda}$, there exists only one $t$ such that $\gamma\in  \hat{\Lambda}_t$. For any $0\leq s\leq t$, the value of
$\gamma_t$ at time $s$ will be denoted by $\gamma_t(s)$. Moreover, if
a path $\gamma_t$ is fixed, the path $\gamma_t|_{[0,s]}$, for $0\leq s\leq t$, will denote the restriction of the path  $\gamma_t$ to the interval
$[0,s]$.
\par
        For convenience, define for $ x\in H,\gamma_t,{\bar{\gamma}}_{\bar{t}}\in \hat{\Lambda}$, $0\leq t\leq \bar{t}\leq T$,
\begin{eqnarray*}
  \gamma^x_{t}(s)&:=&\gamma_t(s){\mathbf{1}}_{[0,t)}(s)+(\gamma_t(t)+x){\mathbf{1}}_{\{t\}}(s), \ \ s\in [0, {t}];\\
  \gamma_{t,\bar{t}}(s)&:=&\gamma_t(s){\mathbf{1}}_{[0,t)}(s)+\gamma_t(t){\mathbf{1}}_{[t,\bar{t}]}(s), \ \ s\in [0,\bar{t}];\\
  \gamma_{t,\bar{t},A}(s)&:=&\gamma_t(s){\mathbf{1}}_{[0,t)}(s)+e^{A(s-t)}\gamma_t(t){\mathbf{1}}_{[t,\bar{t}]}(s), \ \ s\in [0,\bar{t}].
\end{eqnarray*}
\par
 We define a norm   and a metric on $\hat{\Lambda}$ as follows: for any $0\leq t\leq \bar{t}\leq T$ and $\gamma_t,\bar{\gamma}_{\bar{t}}\in \hat{\Lambda}$,
\begin{eqnarray}\label{2.1}
   ||\gamma_t||_0:=\sup_{0\leq s\leq t}|\gamma_t(s)|,\ \ \  \ d_{\infty}(\gamma_t,\bar{\gamma}_{\bar{t}}):=|t-\bar{t}|+||\gamma_{t,T,A}-\bar{\gamma}_{\bar{t},T,A}||_0.
               %+\sup_{0\leq s\leq T}|e^{((s-t)\vee0)A}\gamma_{t}(s\wedge t)-e^{((s-\bar{t})\vee0)A}\bar{\gamma}_{\bar{t}}(s\wedge \bar{t})|.
\end{eqnarray}
%Here and in the sequel,
% for notational simplicity,
%we use $||\gamma_{t}-\bar{\gamma}_{\bar{t}}||_0$ to denote $||\gamma_{t,\bar{t}}-\bar{\gamma}_{\bar{t}}||_0$.
Then $(\hat{\Lambda}_t, ||\cdot||_0)$ is a Banach space, and $(\hat{\Lambda}^t, d_{\infty})$ is a complete metric space by Lemma \ref{lemma2.1111}.
Following Dupire \cite{dupire1}, we define spatial derivatives of $f:\hat{\Lambda}\rightarrow R$, if exist, in the standard sense: if there exists a $B\in H$ such that
$$
\lim_{|h|\rightarrow0}\frac{|f(\gamma_t^{h})-f(\gamma_t)-(B,h)_H|}{|h|}=0,
$$
we say $\partial_{x}f(\gamma_t)=B$,
%\begin{eqnarray}\label{2.2}
% \partial_{x_i}u(\gamma_t):=\lim_{h\rightarrow0}\frac{1}{h}\bigg{[}u(\gamma_t^{he_i})-u(\gamma_t)\bigg{]},\ \
% \partial_{x_ix_j}u:=\partial_{x_i}(\partial_{x_j}u),\ \ i,j=1,2,\cdots,d,
%\end{eqnarray}
and the right time-derivative of $f$, if exists, as:
\begin{eqnarray}\label{2.3}
               \partial_tf(\gamma_t):=\lim_{l\rightarrow0,l>0}\frac{1}{l}\bigg{[}f(\gamma_{t,t+l})-f(\gamma_t)\bigg{]}, \ t<T.
\end{eqnarray}
%For the final time $T$, we define
%$$
%\partial_tu(\gamma_T):=\lim_{t<T,t\uparrow T}\partial_tu(\gamma_t).
%$$
% We take the convention that $\gamma_t$ is column vector, but $\partial_{x}u$ denotes row vector and $\partial_{xx}u$ denotes $d\times d$-matrix.
\begin{definition}\label{definitionc}
       Let $t\in[0,T)$ and $f:\hat{\Lambda}^t\rightarrow R$ be given.
\begin{description}
        \item{(i)}
                 We say $f\in C^0(\hat{\Lambda}^t)$ if $f$ is continuous in $\gamma_s$  on $\hat{\Lambda}^t$ under $d_{\infty}$.
\par
       \item{(ii)}  We say $f\in C^{1}(\hat{\Lambda}^t)\subset C^0(\hat{\Lambda}^t)$ if  $\partial_tf$ and $\partial_{x}f$ exist and are continuous.
\end{description}
\end{definition}
\par
Let $\Lambda_t:= C([0,t],H)$ be the set of all continuous $H$-valued functions defined over $[0,t]$. We denote ${\Lambda}^t=\bigcup_{s\in[t,T]}{\Lambda}_{s}$  and let  ${\Lambda}$ denote ${\Lambda}^0$.
 Clearly, $\Lambda:=\bigcup_{t\in[0,T]}{\Lambda}_{t}\subset\hat{\Lambda}$, and each $\gamma\in \Lambda$ can also be viewed as an element of $\hat{\Lambda}$. $(\Lambda_t, ||\cdot||_0)$ is a
 Banach space, and $(\Lambda^t,d_{\infty})$ is a complete metric space by Lemma \ref{lemma2.1111}.
 $f:\Lambda\rightarrow R$ and $\hat{f}:\hat{\Lambda}\rightarrow R$ are called consistent
  on $\Lambda$ if $f$ is the restriction of $\hat{f}$ on $\Lambda$.
\begin{definition}\label{definitionc2}
       Let  $t\in [0,T)$ and  $f:{\Lambda}^t\rightarrow R$  be given.
\begin{description}
        \item{(i)}
                 We say $f\in C^0({\Lambda}^t)$ if $f$ is continuous in $\gamma_s$  on $\Lambda^t$ under $d_{\infty}$. %For simplicity, we denote $C_0^0({\Lambda})$ by $C^0({\Lambda})$.
\par
       \item{(ii)} We say $f\in C^{1}({\Lambda}^t)$ if
        there exists $\hat{f}\in C^{1}(\hat{{\Lambda}}^t)$ which is consistent with $f$ on $\Lambda^t$.
\end{description}
\end{definition}
 \par
 Let  $(U,d)$ is a metric space. An admissible control  $u=\{u(r),  r\in [t,s]\}$ on $[t,s]$ (with $0\leq t\leq s\leq T$) is a  measurable function taking values in $U$. The set of all admissible controls on $[t,s]$ is denoted by ${\cal{U}}[t,s]$,
      i.e.,
$$
                             {\cal{U}}[t,s]:=\{u(\cdot):[t,s]\rightarrow U|\ u(\cdot) \ \mbox{is
                             measurable}\}.
$$
 Now, we describe some continuous properties of the solutions of  state equation
(\ref{state1}) and value functional (\ref{value1}). First let us assume that functionals
 $ F:{\Lambda}\times U\rightarrow H$,  $
        q: {\Lambda}\times U\rightarrow R$ and $\phi: {\Lambda}_T\rightarrow R
$ satisfy the
following assumption.
      % For any $t\in [0,T]$, denote by ${\cal{H}}^2(t,T)$ the space of all ${{\cal{F}}}^t$-adapted, $R^d$-valued  processes $(Y(s))_{t\leq s\leq T}$ such that
%       $||Y||^2=E[\int^{T}_{t}|Y(s)|^2ds]<\infty$ and by ${\cal{S}}^2(t,T)$ the space of all ${{\cal{F}}}^t$-adapted, $R$-valued continuous processes
%       $(Y(s))_{t\leq s\leq T}$ such that
%       $||Y||^2=E[\sup_{t\leq s\leq T}|Y(s)|^2]<\infty$.
\begin{hyp}\label{hypstate}
\begin{description}
     \item{(i)}
                 The operator $A$ is the generator of a $C_0$ contraction
                 semigroup $\{e^{tA}, t\geq0\}$ of bounded linear operators in the
                 Hilbert space $H$.
                 \par
        \item{(ii)}
        For every fixed $\gamma_t\in\Lambda$, $F(\gamma_t,\cdot)$ and
        $q(\gamma_t,\cdot)$ are continuous in $u$.
\par
       \item{(iii)}
                 There exists a constant $L>0$ such that, for all $(t,\gamma_t,\zeta_T,u)$,  $ (s, \eta_s,\zeta'_Tu) \in [0,T]\times {\Lambda}\times{\Lambda_T}\times U$,
      \begin{eqnarray}\label{assume1111}
                |F(\gamma_t,u)|^2\leq
                 L^2(1+||\gamma_t||_0^2),\ \ \
                 |F(\gamma_t,u)-F(\eta_s,u)|\leq
                 Ld_\infty(\gamma_t,\eta_s);
             \end{eqnarray}
              \begin{eqnarray}\label{asumme2222}
               |q(\gamma_t,u)-q(\eta_s,u)|\leq Ld_\infty(\gamma_t,\eta_s), \ \ |q(\gamma_t,u)|\leq L(1+||\gamma_t||_0);
\end{eqnarray}
 \begin{eqnarray}\label{asumme3333}
                 |\phi(\zeta_T)-\phi(\zeta'_T)|\leq L||\zeta-\zeta'_T||_0, \ \ \ |\phi(\zeta_T)|\leq L(1+||\zeta||_0).
\end{eqnarray}
\end{description}
\end{hyp}
\par
              We say that $X$ is a mild solution of equation $(\ref{state1})$ if $X\in C^0(\Lambda)$ and it satisfies:
\begin{eqnarray*}
            X(s)=e^{(s-t)A}\gamma_t(t)+\int_{t}^{s}{e^{(s-\sigma)A}}F(X_\sigma,u(\sigma))d\sigma,  \ s\in [t,T];\ \ \mbox{and} \ \ X(s)=\gamma_t(s),  \ s\in[0,t).
\end{eqnarray*}
\par
The following lemma is standard.
\begin{lemma}\label{lemmaexist}
\ \ Assume that Hypothesis \ref{hypstate}  (iii)  holds. Then for every $u\in {\cal{U}}[t,T]$,
$\gamma_t\in {\Lambda}$, (\ref{state1}) admits a
unique mild solution $X^{\gamma_t,u}$.  Moreover, if we let  $X^{\eta_t,u}$  be the solutions of  (\ref{state1})
 corresponding $\eta_t\in \Lambda$ and $u\in {\cal{U}}[t,T]$. Then the following estimates hold:
\begin{eqnarray}\label{state1est}
               ||X_T^{\gamma_t,u}-X_T^{\eta_t,u}||_0\leq C_1||\gamma_t-\eta_t||_0,\ \ \ \
                ||X_T^{\gamma_t,u}||_0\leq C_1(1+||\gamma_t||_0).
\end{eqnarray}
              The constant $C_1$ depending only on   $T$, $L$ and
              $M_1=:\sup_{s\in [0,T]}|e^{sA}|$.
\end{lemma}
{\bf  Proof  }. \ \
By Picard iteration, we can obtain the existence and the uniqueness
of the mild solution. By Gronwall's inequality, together with assumption
(\ref{assume1111}), we can prove (\ref{state1est}). \ \ $\Box$
  \par
  The next result contains the local boundedness and the continuity of the trajectory $X^{\gamma_t,u}$ and  value functional $V$. In what follows, $C$ is an absolute
constant, that can be different in different places.
  \begin{lemma}\label{lemmaexist111}
\ \ Assume that Hypothesis \ref{hypstate}  (iii)  holds. Then, for any $0\leq t\leq \bar{t}\leq T$, $\gamma_t, \eta_t\in {\Lambda}$ and $u\in {\cal{U}}[t,T]$,
\begin{eqnarray}\label{2.6}
            \sup_{u\in {\cal{U}}[t,T]}|X^{\gamma_t,u}(s)-e^{(s-t)A}\gamma_t(t)|\leq
            C(1+||\gamma_t||_0)|s-t|, \ \ \ s\in[t,T];
\end{eqnarray}
\begin{eqnarray}\label{0604}
||X^{\eta_t,u}_T-X^{\gamma_{t,\bar{t},A},u}_T||_0\leq C(1+||\eta_t||_0)(\bar{t}-t)+C||\eta_t-\gamma_t||_0;
\end{eqnarray}
\begin{eqnarray}\label{3.5}
 |V(\gamma_t)|\leq C(1+||\gamma_t||_0);
 \end{eqnarray}
\begin{eqnarray}\label{0604001}
|V(\gamma_{t,\bar{t},A})-V(\eta_t)|
                        \leq C(1+||\eta_t||_0)(\bar{t}-t)+C||\eta_t-\gamma_t||_0.
\end{eqnarray}
\end{lemma}
{\bf  Proof  }. \ \
 For any  $\gamma_t\in\Lambda$, by (\ref{assume1111}) and (\ref{state1est}), we
                obtain the following result:
\begin{eqnarray*}
                   |X^{\gamma_t,u}(s)-e^{(s-t)A}\gamma_t(t)|
               \leq LM_1(1+C_1(1+||\gamma_t||_0))|s-t|.
\end{eqnarray*}
    Taking the supremum in ${\cal{U}}[t,T]$, we obtain (\ref{2.6}).
  For any $0\leq t\leq \bar{t}\leq T$, $\gamma_t, \eta_t\in {\Lambda}$ and $u\in {\cal{U}}[t,T]$, by (\ref{assume1111}) and (\ref{state1est}), we have
%\begin{eqnarray*}
%             |X^{\gamma_t,u}(s)-e^{A(s-t)}\gamma_t(t))|\leq \int^{s}_{t}L(1+||X^{\gamma_t,u}_r||_0)dr\leq C(1+||\gamma_t||_0)(s-t).
%\end{eqnarray*}
\begin{eqnarray*}
             &&\sup_{\bar{t}\leq s\leq \sigma}|X^{\eta_t,u}(s)-X^{\gamma_{t,\bar{t},A},u}(s)|\\
             &\leq& M_1|\eta_t(t)-\gamma_t(t)|+%\sup_{\bar{t}\leq s\leq \sigma}
             \int^{\bar{t}}_{t}|e^{(s-r)A}F(X^{\eta_t,u}_r,u(r))|dr\\
             &&
                                                             +\sup_{\bar{t}\leq s\leq \sigma}\int^{s}_{\bar{t}}|e^{(s-r)A}(F(X^{\eta_t,u}_r,u(r)))-F(X^{\gamma_{t,\bar{t},A},u}_r,u(r))|dr\\
                                                            % &\leq& M|\eta_t(t)-\gamma_t(t)|+LM\int^{\bar{t}}_{t}(1+||X^{\eta_t,u}_r||_0)dr
%                                                             +LM\int^{\sigma}_{\bar{t}}||X^{\eta_t,u}_r-X^{\gamma_{t,\bar{t},A},u}_r||_0dr\\
                                                              &\leq&M_1|\eta_t(t)-\gamma_t(t)|+LM_1(1+C_1(1+||\eta_t||_0))(\bar{t}-t)
                                                             +LM_1\int^{\sigma}_{\bar{t}}||X^{\eta_t,u}_r-X^{\gamma_{t,\bar{t},A},u}_r||_0dr.
\end{eqnarray*}
Thus,
%\begin{eqnarray*}
%             ||X^{\eta_t,u}_\sigma-X^{\gamma_{t,\bar{t},A},u}_\sigma||_0
%                                                              \leq M_1||\eta_t-\gamma_t||_0+LM_1(1+C_1(1+||\eta_t||_0))(\bar{t}-t)
%                                                             +LM_1\int^{\sigma}_{\bar{t}}||X^{\eta_t,u}_r-X^{\gamma_{t,\bar{t},A},u}_r||_0dr.
%\end{eqnarray*}
\begin{eqnarray*}
             ||X^{\eta_t,u}_\sigma-X^{\gamma_{t,\bar{t},A},u}_\sigma||_0
                                                              \leq C||\eta_t-\gamma_t||_0+C(1+||\eta_t||_0)(\bar{t}-t)
                                                             +C\int^{\sigma}_{\bar{t}}||X^{\eta_t,u}_r-X^{\gamma_{t,\bar{t},A},u}_r||_0dr.
\end{eqnarray*}
Then, by Gronwall's inequality, we obtain (\ref{0604}). Next, by (\ref{asumme2222}), (\ref{asumme3333}) and (\ref{0604}), we get
\begin{eqnarray*}
                        &&|J(\gamma_{t,\bar{t},A},u)-J(\eta_t,u)|\\&\leq& L\int_{t}^{\bar{t}}(1+||X^{\eta_t,u}_\sigma||_0)d\sigma
                        +L\int_{\bar{t}}^{T}||X^{\gamma_{t,\bar{t},A},u}_\sigma-X^{\eta_t,u}_\sigma||_0d\sigma+L||X^{\gamma_{t,\bar{t},A},u}_T-X^{\eta_t,u}_T||_0\\
                        &\leq& C(1+||\eta_t||_0)(\bar{t}-t)+C||\eta_t-\gamma_t||_0.
\end{eqnarray*}
Thus, taking the infimum in $u(\cdot)\in {\cal{U}}[t,T]$, we obtain (\ref{0604001}). By the similar procedure, we can show
                     (\ref{3.5}) holds true. The lemma  is proved. \ \  $\Box$
\par
  Next, we present the dynamic programming
                        principle (DPP) for optimal control problems (\ref{state1}) and (\ref{value1}).
\begin{theorem}\label{theorem3.3}
                      Assume the Hypothesis \ref{hypstate} (ii) and (iii) hold true. Then, for every $(t,\gamma_t)\in [0,T)\times \Lambda$ and $s\in [t,T]$, we have that
\begin{eqnarray}\label{3.7}
             V(\gamma_t)=\inf_{u\in {\cal{U}}[t,T]}\bigg{[}\int_{t}^{s}q(X^{\gamma_t,u}_\sigma,u(\sigma))d\sigma+V(X^{\gamma_t,u}_s)\bigg{]}.
\end{eqnarray}
\end{theorem}
The proof is very similar to the case without path-dependent (see Theorem 1.1  of Chapter 6 in page 224 of \cite{li}). For the convenience of readers, here we give its proof.
 \par
{\bf  Proof}. \ \ First of all, for any  $u\in{\cal{U}}[s,T]$,
$s\in[t,T]$ and any  $u\in{\cal{U}}[t,s]$, by putting them
             concatenatively, we get $u\in{\cal{U}}[t,T]$.
                  Let us denote the right-hand side of (\ref{3.7}) by
                  $\overline{V}({\gamma_t})$. By (\ref{value1}), %we have
$$
           V({\gamma_t})\leq J({\gamma_t,u})
              =\int_{t}^{s}q(X^{\gamma_t,u}_\sigma,u(\sigma))d\sigma+J(X^{{\gamma_t,u}}_s,{{u}}),\
              u(\cdot)\in {\cal{U}}[t,T].
$$
             Thus, taking the infumum over $u(\cdot)\in
             {\cal{U}}[s,T]$, we obtain
$$
           V({\gamma_t})\leq \int_{t}^{s}q(X^{\gamma_t,u}_\sigma,u(\sigma))d\sigma+V(X^{{\gamma_t,u}}_s).
$$
             Consequently,
$$
                 V({\gamma_t})\leq \overline{V}({\gamma_t}).
$$
            On the other hand, for any $\varepsilon>0$, there exists a $u^\varepsilon\in {\cal{U}}[t,T]$ such that
       \begin{eqnarray*}
           V({\gamma_t})+\varepsilon&\geq& J({\gamma_t},u^\varepsilon)
           =
              \int_{t}^{s}q(X^{\gamma_t,u^\varepsilon}_\sigma,{u^\varepsilon}(\sigma))d\sigma
                           +J(X^{\gamma_t,u^\varepsilon}_s,{u^\varepsilon})\\
              &\geq&\int_{t}^{s}q(X^{\gamma_t,u^\varepsilon}_\sigma,{u^\varepsilon}(\sigma))d\sigma
                           +V(X^{\gamma_t,u^\varepsilon}_s)\geq\overline{V}({\gamma_t}).
\end{eqnarray*}
             Hence, (\ref{3.7}) follows.\ \ $\Box$
             \par
                The following
                lemma is needed to prove the existence  of viscosity solutions.
\begin{theorem}\label{theoremito}
\ \
Suppose $X$ is a solution of (\ref{state1}), $\varphi\in C^{1}({\Lambda}^{{\bar{t}}})$  and $A^*\partial_x\varphi\in C^0({\Lambda}^{{\bar{t}}})$ for some $\bar{t}\in[t,T)$. Then for any $s\in [\bar{t},T]$:
\begin{eqnarray}\label{statesop0}
                 \varphi(X_s)=\varphi(X_{\bar{t}})+\int_{\bar{t}}^{s}\partial_t\varphi(X_\sigma)d\sigma+\int^{s}_{\bar{t}}(A^*\partial_x\varphi(X_\sigma),X(\sigma))_H+(\partial_x\varphi(X_\sigma),F(X_\sigma,u(\sigma))_Hd\sigma.
\end{eqnarray}
\end{theorem}
\par
{\bf  Proof}. \ \  Denote $X^n=X{\mathbf{1}}_{[0,\bar{t})}+\sum^{2^n-1}_{i=0}X(t_{i+1}){\mathbf{1}}_{[t_i,t_{i+1})}+X(s){\mathbf{1}}_{\{s\}}$ which is a c$\grave{a}$dl$\grave{a}$g piecewise constant approximation of $X$.
Here $t_i=\bar{t}+\frac{i(s-\bar{t})}{2^n}$. For every $\gamma_\sigma\in \hat{\Lambda}$, define $\gamma_{\sigma-}\in \hat{\Lambda}$ by
$$
               \gamma_{\sigma-}(\theta)=\gamma_{\sigma}(\theta),\ \ \theta\in [0,\sigma),\ \ \mbox{and}\  \ \gamma_{\sigma-}(\sigma)=\lim_{\theta\uparrow \sigma}\gamma_{\sigma}(\theta).
$$ We start with the decomposition
\begin{eqnarray}\label{decom}
                   \varphi({X^n_{t_{i+1}}}_{-})-\varphi({X^n_{t_{i}}}_{-})=\varphi({X^n_{t_{i+1}}}_{-})-\varphi({X^n_{t_{i}}})
                   +\varphi({X^n_{t_{i}}})-\varphi({X^n_{t_{i}}}_{-}).
\end{eqnarray}
 Let $\psi(l)=\varphi({X^n_{t_{i},t_i+l}})$, we have $\varphi({X^n_{t_{i+1}}}_{-})-\varphi({X^n_{t_{i}}})=\psi(h)-\psi(0)$, where $h=\frac{s-\bar{t}}{2^n}$.
 Since $\varphi\in C^{1}({\Lambda}^{\bar{t}})$,  the right derivative of $\psi$ is continuous, therefore,
 $$
                                \varphi({X^n_{t_{i+1}}}_{-})-\varphi({X^n_{t_{i}}})=\int^{t_{i+1}}_{t_i}\partial_t\varphi(X_{t_i,l}^n)dl.
$$
The term $\varphi({X^n_{t_{i}}})-\varphi({X^n_{t_{i}}}_{-})$ in (\ref{decom}) can be written $\pi(X(t_{i+1}))-\pi(X(t_i))$, where
$\pi(l)=\varphi({X^n_{t_{i}}}_{-}+(l-X(t_i)){\mathbf{1}}_{\{t_{i}\}})$. Since $\varphi\in C^{1}({\Lambda}^{\bar{t}})$, $\pi$ is a $C^1$ function and $\nabla_x\pi(l)=\partial_x\varphi({X^n_{t_{i}}}_{-}+(l-X(t_i)){\mathbf{1}}_{\{t_{i}\}})$.
Thus, by Proposition 5.5 in Chapter 2 of \cite{li}, we have that:
\begin{eqnarray*}
                              \pi({X(t_{i+1})})-\pi(X(t_i))&=&\int^{t_{i+1}}_{t_i}(A^*\partial_x\varphi({X^n_{t_{i}}}_{-}+(X(\sigma)-X(t_i)){\mathbf{1}}_{\{t_{i}\}}),X(\sigma))_H\\
                              &&~~~~~+(\partial_x\varphi({X^n_{t_{i}}}_{-}+(X(\sigma)-X(t_i)){\mathbf{1}}_{\{t_{i}\}}),F(X_\sigma,u(\sigma)))_Hd\sigma.
\end{eqnarray*}
 Summing over $i\geq 0$ and denoting $i(\sigma)$ the index such that $\sigma\in [t_{i(\sigma)},t_{i(\sigma)+1})$, we obtain
\begin{eqnarray*}
                                 \varphi(X^n_s)-\varphi(X^n_{\bar{t}})
                                &=&\int^{s}_{\bar{t}}\partial_t\varphi(X_{t_{i(\sigma),\sigma}}^n)d\sigma+\int^{s}_{\bar{t}}(A^*\partial_x\varphi({X^n_{t_{i(\sigma)}}}_{-}+(X(\sigma)-X(t_{i(\sigma)})){\mathbf{1}}_{\{t_{i(\sigma)}\}}),X(\sigma))_H\\
                              &&~~~~~+(\partial_x\varphi({X^n_{t_{i(\sigma)}}}_{-}+(X(\sigma)-X(t_{i(\sigma)})){\mathbf{1}}_{\{t_{i(\sigma)}\}}),F(X_\sigma,u(\sigma)))_Hd\sigma.
\end{eqnarray*}
  $\varphi(X^n_s)$ and $\varphi(X^n_{\bar{t}})$ converge to $\varphi(X_s)$ and $\varphi(X_{\bar{t}})$, respectively. Since all approximations of $X$ appearing in the  integrals have a $||\cdot||_{0}$-distance from $X_s$ less than $||X^n_s-X_s||_0\rightarrow0$, $\varphi\in C^{1}({\Lambda}^{{t}})$ and  $A^*\partial_x\varphi\in C^0({\Lambda}^{{\bar{t}}})$ imply that the integrands appearing in the above integrals converge respectively to $\partial_t\varphi(X_\sigma)$, $A^*\partial_x\varphi(X_\sigma)$ and  $\partial_x\varphi(X_\sigma)$ as $n\rightarrow\infty$. By $X$ is continuous, and $\varphi\in C^{1}({\Lambda}^{\bar{t}})$ and  $A^*\partial_x\varphi\in C^0({\Lambda}^{{\bar{t}}})$, the integrands in the various above integrals are bounded.  The dominated convergence then ensure that the Lebesgue integrals converge to the terms appearing in (\ref{statesop0}) as $n\rightarrow\infty$.\ \ $\Box$
\par
We conclude this section with  the following four    lemmas which will be used to prove the uniqueness and  stability  of viscosity solutions.
\begin{definition}\label{gaupe}
              Let $t\in [0,T]$ be fixed.  We say that a continuous functional $\rho:\Lambda^t\times \Lambda^t\rightarrow [0,+\infty)$ is a {gauge-type function} provided that:
             \begin{description}
        \item{(i)} $\rho(\gamma_s,\gamma_s)=0$ for all $(s,\gamma_s)\in [t,T]\times \Lambda^t$,
        \item{(ii)} for any $\varepsilon>0$, there exists $\delta>0$ such that, for all $\gamma_s, \eta_l\in \Lambda^t$, we have $\rho(\gamma_s,\eta_l)\leq \delta$ implies that
        $d_\infty(\gamma_s,\eta_l)<\varepsilon$.
        \end{description}
\end{definition}
The following lemma is a modification of Borwein-Preiss variational principle (see Theorem 2.5.2 in  Borwein \& Zhu  \cite{bor1}). It will be used to
get a maximum of a perturbation of the auxiliary function in the proof of uniqueness. The proof  is completely similar to the finite dimensional case (see Lemma 2.12 in \cite{zhou5}). Here we omit it.
\begin{lemma}\label{theoremleft} %(Borwein-Preiss Variational Principle)
Let $t\in [0,T]$ be fixed and let $f:\Lambda^t\rightarrow R$ be an upper semicontinuous functional  bounded from above. Suppose that $\rho$ is a gauge-type function
 and $\{\delta_i\}_{i\geq0}$ is a sequence of positive number, and suppose that $\varepsilon>0$ and $(t_0,\gamma^0_{t_0})\in [t,T]\times \Lambda^t$ satisfy
 $$
f(\gamma^0_{t_0})\geq \sup_{(s,\gamma_s)\in [t,T]\times \Lambda^t}f(\gamma_s)-\varepsilon.
 $$
 Then there exist $(\hat{t},\hat{\gamma}_{\hat{t}})\in [t,T]\times \Lambda^t$ and a sequence $\{(t_i,\gamma^i_{t_i})\}_{i\geq1}\subset [t,T]\times \Lambda^t$ such that
  \begin{description}
        \item{(i)} $\rho(\gamma^0_{t_0},\hat{\gamma}_{\hat{t}})\leq \frac{\varepsilon}{\delta_0}$,  $\rho(\gamma^i_{t_i},\hat{\gamma}_{\hat{t}})\leq \frac{\varepsilon}{2^i\delta_0}$ and $t_i\uparrow \hat{t}$ as $i\rightarrow\infty$,
        \item{(ii)}  $f(\hat{\gamma}_{\hat{t}})-\sum_{i=0}^{\infty}\delta_i\rho(\gamma^i_{t_i},\hat{\gamma}_{\hat{t}})\geq f(\gamma^0_{t_0})$, and
        \item{(iii)}  $f(\gamma_s)-\sum_{i=0}^{\infty}\delta_i\rho(\gamma^i_{t_i},\gamma_s)
            <f(\hat{\gamma}_{\hat{t}})-\sum_{i=0}^{\infty}\delta_i\rho(\gamma^i_{t_i},\hat{\gamma}_{\hat{t}})$ for all $(s,\gamma_s)\in [\hat{t},T]\times \Lambda^{\hat{t}}\setminus \{(\hat{t},\hat{\gamma}_{\hat{t}})\}$.

        \end{description}
\end{lemma}
  \par
   We define, for every $\gamma_t\in \hat{\Lambda}$,
 \begin{eqnarray*}
S(\gamma_t)=\begin{cases}
            \frac{(||\gamma_{t}||_0^2-|\gamma_{t}(t)|^2)^2}{||\gamma_{t}||^2_0}, \
         ~~ ||\gamma_{t}||_0\neq0; \\
0, \ ~~~~~~~~~~~~~~~~~~~ ||\gamma_{t}||_0=0.
\end{cases}
\end{eqnarray*}
\begin{lemma}\label{theoremS}
 $S(\cdot)\in C^{1}(\hat{\Lambda})$. Moreover,
\begin{eqnarray}\label{s0}
                ||\gamma_t||_0^2\leq   S(\gamma_t)+   2|\gamma_t(t)|^2\leq 3||\gamma_t||_0^2.
\end{eqnarray}
\end{lemma}
\par
   {\bf  Proof  }. \ \  First, we prove $S\in C^0(\hat{\Lambda})$. %Let $\eta_s\rightarrow \gamma_t$ in $(\Lambda,d_\infty)$.
   For any $\gamma_t,\eta_s\in \hat{\Lambda}$, if $s\geq t$,
   $$
   |\gamma_t(t)-\eta_s(s)|\leq |\gamma_t(t)-e^{(s-t)A}\gamma_t(t)|+|e^{(s-t)A}\gamma_t(t)-\eta_s(s)|, %\rightarrow0, \ \ \mbox{as}\ d_\infty(\gamma_t,\eta_s)\rightarrow0,
   $$
   and
   \begin{eqnarray*}
   |||\gamma_t||_0-||\eta_s||_0|&\leq& ||\gamma_{t,s,A}||_0-||\gamma_t||_0+||\gamma_{t,s,A}-\eta_s||_0\\
   &\leq& \sup_{t\leq l\leq s}|(e^{(l-t)A}-I)\gamma_t(t)|+||\gamma_{t,s,A}-\eta_s||_0;
   %\rightarrow0 \ \ \mbox{as}\ d_\infty(\gamma_t,\eta_s)\rightarrow0.
\end{eqnarray*}
   %Thus we have $S(\eta_s)\rightarrow S(\gamma_t)$ as $\ d_\infty(\gamma_t,\eta_s)\rightarrow0$.
    %Otherwise, we may assume $s<t$, then
    if $s<t$,
   \begin{eqnarray*}
   |\gamma_t(t)-\eta_s(s)|&\leq& |\gamma_t(t)-e^{(t-s)A}\eta_s(s)|+|e^{(t-s)A}(\eta_s(s)-\gamma_t(s))|+|e^{(t-s)A}(\gamma_t(s)-\gamma_t(t-))|\\
   &&+|(e^{(t-s)A}-I)\gamma_t(t-)|+|\gamma_t(t-)-\gamma_t(s)|+|\gamma_t(s)-\eta_s(s)|, %\rightarrow0, \ \ \mbox{as}\ d_\infty(\gamma_t,\eta_s)\rightarrow0,
 \end{eqnarray*}
 and
 \begin{eqnarray*}
   &&|||\gamma_t||_0-||\eta_s||_0|\leq ||\eta_{s,t,A}||_0-||\eta_s||_0+||\eta_{s,t,A}-\gamma_t||_0\\
   &\leq& \sup_{s\leq l\leq t}|(e^{(l-s)A}-I)\gamma_t(t)|+\sup_{s\leq l\leq t}|(e^{(l-s)A}-I)||\gamma_t(t)-\eta_s(s)|+||\eta_{s,t,A}-\gamma_t||_0.
   %\rightarrow0 \ \ \mbox{as}\ d_\infty(\gamma_t,\eta_s)\rightarrow0.
\end{eqnarray*}
    Then we have $S(\eta_s)\rightarrow S(\gamma_t)$ as $\eta_s\rightarrow\gamma_t$ under $d_\infty$. Thus $S\in C^0(\hat{\Lambda})$.
   Second, by the definition of $S(\cdot)$, it is clear that $\partial_tS(\gamma_{t})=0$ for all $\gamma_t\in \hat{\Lambda}$.
  Next, we consider $ \partial_{x}S$. %For every $\gamma_t\in \hat{\Lambda}$, if $||\gamma_t||_0=0$,
%  \begin{eqnarray}\label{ss4}
% 0\leq\lim_{h\rightarrow0}\bigg{|}\frac{S(\gamma_t^{h})-S(\gamma_t)}{h}\bigg{|}
%  %\leq  \lim_{h\rightarrow0}\bigg{|}\frac{{(|\gamma_t(t)|^2-|\gamma_t(t)+{he_i}|^2)^3}}{h(1+|\gamma_t|_C^4)} \bigg{|}
%  \leq\lim_{h\rightarrow0}\bigg{|}\frac{h^4}{h^3} \bigg{|}=0;
%   \end{eqnarray}
%
%  If $||\gamma_t||_0\neq0$, for every $h\in H$,
%   \begin{eqnarray*}
%   S(\gamma_t^{h})-S(\gamma_{t})
%   =\frac{(||\gamma^{h}_{t}||_0^2-|\gamma_t(t)+{h}|^2)^2}{||\gamma^{h}_{t}||_0^2}
%   -\frac{(||\gamma_{t}||_0^2-|\gamma_t(t))|^2)^2}{||\gamma_{t}||_0^2}.
%   \end{eqnarray*}
For every $\gamma_t\in \hat{\Lambda}$, let $||\gamma_t||_{0^-}=\sup_{0\leq s<t}|\gamma_t(s)|$. %  and $\gamma^i_t(t)=\gamma_t(t)e_i,\ i=1,2,\ldots, d$.
Then,
 if $|\gamma_t(t)|<||\gamma_t||_{0^-}$,
\begin{eqnarray*}
  &&\lim_{|h|\rightarrow0}\frac{\bigg{|}S(\gamma_t^{h})-S(\gamma_{t})+\frac{4(||\gamma_t||^2_{0}-|\gamma_t(t)|^2)(\gamma_t(t),h)_H}{||\gamma_t||^2_{0}}\bigg{|}}{|h|}\\
  &=&\lim_{|h|\rightarrow0}\frac{{|(||\gamma_t||^2_{0}-|\gamma_t(t)+{h}|^2)^2}
   -{(||\gamma_t||^2_{0}-|\gamma_t(t)|^2)^2}+4(||\gamma_t||^2_{0}-|\gamma_t(t)|^2)(\gamma_t(t),h)_H|}{|h|\times||\gamma_t||^2_{0}}\nonumber\\
   &=&\lim_{|h|\rightarrow0}\frac{{|-(2||\gamma_t||^2_{0}-|\gamma_t(t)+{h}|^2-|\gamma_t(t)|^2)(2(\gamma_t(t),h)_H+|h|^2)+4(||\gamma_t||^2_{0}-|\gamma_t(t)|^2)(\gamma_t(t),h)_H|}
  }{|h|\times||\gamma_t||^2_{0}}\nonumber\\
   &=&0.
   \end{eqnarray*}
   Thus,
\begin{eqnarray}\label{s1}
   \partial_xS(\gamma_{t})=-\frac{4(||\gamma_t||^2_{0}-|\gamma_t(t)|^2)
    \gamma_t(t)}{||\gamma_t||^2_{0}}.
  \end{eqnarray}
 If  $|\gamma_t(t)|>||\gamma_t||_{0^-}$,
\begin{eqnarray}\label{s2}
   \partial_{x}S(\gamma_{t})=0;
   \end{eqnarray}
if  $|\gamma_t(t)|=||\gamma_t||_{0^-}\neq0$,
since
\begin{eqnarray*}
||\gamma^{h}_t||_0^2-|\gamma_t(t)+{h}|^2=
\begin{cases}
0,\ \ \ \ \ \ \ \ \ \  \ \ \ \ \ \ \ \ \ \ \  \ \ \ \ \ \  |\gamma_t(t)+h|\geq |\gamma_t(t)|,\\
  |\gamma_t(t)|^2-|\gamma_t(t)+{h}|^2,\ \ |\gamma_t(t)+h|<|\gamma_t(t)|,
\end{cases}
\end{eqnarray*}
we have%, %if $|\gamma_t(t)-a_{\hat{t}}(\hat{t})|=||\gamma_t-a_{\hat{t},t}||_{0^-}\neq0$,
\begin{eqnarray}\label{s3}
 0\leq\lim_{|h|\rightarrow0}\frac{|S(\gamma_t^{h})-S(\gamma_t)|}{|h|}
  %\leq  \lim_{h\rightarrow0}\bigg{|}\frac{{(|\gamma_t(t)|^2-|\gamma_t(t)+{he_i}|^2)^3}}{h(1+|\gamma_t|_C^4)} \bigg{|}
  \leq\lim_{|h|\rightarrow0}\frac{|h|^2{(|h|+2|\gamma_t(t)|)^2}}{|h|||\gamma^{h}_t||_0^2} =0;
   \end{eqnarray}
    if $|\gamma_t(t)-a_{\hat{t}}(\hat{t})|=||\gamma_t-a_{\hat{t},t}||_{0^-}=0$,
\begin{eqnarray}\label{ss4}
 \partial_{x}S(\gamma_{t})=0.
   \end{eqnarray}
From (\ref{s1}), (\ref{s2}), (\ref{s3}) and (\ref{ss4}) we obtain that
\begin{eqnarray*}
    \partial_{x}S(\gamma_t)=\begin{cases}-\frac{4(||\gamma_t||^2_{0}-|\gamma_t(t)|^2)
    \gamma_t(t)}{||\gamma_t||^2_{0}}, \ \ \ \  ||\gamma_t||^2_{0}\neq0,\\
    0, ~~~~~~~~~~~~~~~~~~~~~~~~~~~~||\gamma_t||^2_{0}=0.
    \end{cases}
\end{eqnarray*}
It is clear that $\partial_{x}S\in C^0(\hat{\Lambda})$. Thus, we have show that $S(\cdot)\in C^{1}(\hat{\Lambda})$.
%\par
%For every $g\in C^1(R)$, there exists $0<\tau<1$ such that
%$$
%g(S(\gamma^h_t))-g(S(\gamma_t))=g'(S(\gamma_t)+\tau(S(\gamma^h_t)-S(\gamma_t)))(S(\gamma^h_t)-S(\gamma_t)).
%$$
%  Then
%  \begin{eqnarray*}
%  &&\lim_{|h|\rightarrow0}\frac{|g(S(\gamma^h_t))-g(S(\gamma_t))-g'(S(\gamma_t))(\partial_{x}S(\gamma_t),h)|}{|h|}\\
%  &\leq&\lim_{|h|\rightarrow0}\frac{|g'(S(\gamma_t)+\tau(S(\gamma^h_t)-S(\gamma_t)))-g'(S(\gamma_t)||(S(\gamma^h_t)-S(\gamma_t))|}{|h|}\\
%  &&+\lim_{|h|\rightarrow0}\frac{|g'(S(\gamma_t)||S(\gamma^h_t)-S(\gamma_t)-(\partial_{x}S(\gamma_t),h)|}{|h|}
%  =0.
%  \end{eqnarray*}
%   Therefore, $$
%   \partial_{x}g(S(\gamma_t))=g'(S(\gamma_t))\partial_{x}S(\gamma_t).
%   $$
%   It is clear that $\partial_{x}g(S(\gamma_t))\in C^0(\Lambda)$. Thus, we have show that $g(S(\cdot))\in C^{1}(\Lambda)$.
\par
Now we prove (\ref{s0}). It is clear that
$$
                      S(\gamma_t)+   2|\gamma_t(t)|^2\leq 3||\gamma_t||_0^2.
$$
On the other hand, if $||\gamma_t||_0\neq0$,
\begin{eqnarray*}
                      S(\gamma_t)+   2|\gamma_t(t)|^2=\frac{(||\gamma_t||_0^2-|\gamma_t(t)|^2)^2-2||\gamma_t||_0^2|\gamma_t(t)|^2}{||\gamma_t||_0^2}
                      =\frac{|\gamma_t(t)|^4+||\gamma_t||_0^4}{||\gamma_t||_0^2}\geq ||\gamma_t||_0^2.
\end{eqnarray*}
% if $||\gamma_t||_0^2=0$, it is clear that
%$$
%                      S(\gamma_t)+   2|\gamma_t(t)|^2\geq||\gamma_t||_0^2.
%$$
Thus, we  have (\ref{s0}) holds true.
 The proof is now complete. \ \ $\Box$
\par
Define, for every $M\in R$,
$$
         \Upsilon^M(\gamma_t):=S(\gamma_t)+M|\gamma_t(t)|^2, \ \  \gamma_t\in \Lambda;
$$
$$
         \Upsilon^M(\gamma_t,\eta_s)= \Upsilon^M(\eta_s,\gamma_t):=\Upsilon^M(\eta_s-\gamma_{t,s,A}), \ \  0\leq t\leq s\leq T, \ \gamma_t,\eta_s\in \Lambda;
$$
and
$$
\overline{\Upsilon}^M(\gamma_t,\eta_s)=\overline{\Upsilon}^M(\eta_s,\gamma_t):=\Upsilon^M(\eta_s,\gamma_t)+|s-t|^2, \ \  0\leq t\leq s\leq T, \ \gamma_t,\eta_s\in \Lambda.
$$
The proof of the following Lemma is completely similar to the finite dimensional case (see Lemma 2.13 in \cite{zhou5}). Here we omit it.
\begin{lemma}\label{theoremS00044} For $M\geq2$, we have that
\begin{eqnarray}\label{jias5}
2\Upsilon^M(\gamma_t)+2\Upsilon^M(\eta_t)\geq\Upsilon^M(\gamma_t+\eta_t), \ \ (t,\gamma_t, \eta_t)\in [0,T]\times\hat{ \Lambda}\times \hat{\Lambda}.
\end{eqnarray}
\end{lemma}

\begin{lemma}\label{theoremS000}  Assume the Hypothesis \ref{hypstate}  holds true. For every $t\in [0,T)$, $\eta_t\in \Lambda$ and $M\geq2$, %and $g\in C^1(R)$ with $g'\geq0$,
we have that
\begin{eqnarray}\label{jias51}
                  \Upsilon^M(X^{\gamma_t,u}_s-\eta_{t,s,A})\leq \Upsilon^M(X^{\gamma_t,u}_t-\eta_t)+\int^{s}_{{t}}(\partial_x\Upsilon^M(X^{\gamma_t,u}_\sigma-\eta_{t,\sigma,A}),F(X^{\gamma_t,u}_\sigma,u(\sigma))_Hd\sigma.
\end{eqnarray}
\end{lemma}
\par
   {\bf  Proof  }. \ \ Let $A_\mu=\mu A(\mu I-A)^{-1}$ be the Yosida approximation of $A$ and
let $X^\mu$ be the solution of the following:
$$
X^\mu(s)=e^{(s-t)A_\mu}\gamma_t(t)+\int^{s}_{t}e^{(s-\sigma)A_\mu}F(X^{\mu}_\sigma,u(\sigma))d\sigma, \ s\in [t,T];\ \  X^\mu(s)=\gamma_t(s), \ s\in[0,t).
$$
 Define $y^\mu$ by $y^\mu(s)=X^\mu(s)-e^{(s-t)A_\mu}\eta_t(t),\ \ t\leq s\leq T$ and $y^\mu(s)=\gamma_t(s)-\eta_t(s),\ \ 0\leq s<t$, then $y^\mu$ be the solution of the following:
$$
                      y^\mu(s)=e^{(s-t)A_\mu}y^\mu_t(t)+\int^{s}_{t}e^{(s-\sigma)A_\mu}F(X^{\mu}_\sigma, u(\sigma))d\sigma,\ s\in [t,T];\ \  y^\mu(s)=\gamma_t(s)-\eta_t(s), \ s\in[0,t).
$$
 By Theorem \ref{theoremito} and Lemma \ref{theoremS}, we have, %for every $g\in C^1(R)$ with $g'\geq0$,
\begin{eqnarray*}
                  \Upsilon^M(y^\mu_s)=\Upsilon^M(y^\mu_{{t}})+\int^{s}_{{t}}(\partial_x\Upsilon^M(y^\mu_\sigma),A_\mu y^\mu(\sigma)+F(X^{\mu}_\sigma,u(\sigma))_Hd\sigma.
\end{eqnarray*}
   Noting that $A$ is the infinitesimal generator of a $C_0$ contraction semigroup, we have, if $||y^\mu_\sigma||^2_{0}\neq0$,
   $$\bigg{(}2My^\mu(\sigma)-\frac{4(||y^\mu_\sigma||^2_{0}-|y^\mu(\sigma)|^2) y^\mu(\sigma)}{||y^\mu_\sigma||^2_{0}},A_\mu y^\mu(\sigma)\bigg{)}_H\leq 0,\ \ \mbox{for}\  M\geq2.$$ Thus,
   \begin{eqnarray*}
                  \Upsilon^M(y^\mu_s)\leq \Upsilon^M(y^\mu_{{t}})+\int^{s}_{{t}}(\partial_x\Upsilon^M(y^\mu_\sigma),F(X^{\mu}_\sigma,u(\sigma))_Hd\sigma.
\end{eqnarray*}
Letting $\mu\rightarrow\infty$, by Proposition 5.4 in Chapter 2 of \cite{li}, we obtain
\begin{eqnarray*}
                  \Upsilon^M(y_s)\leq \Upsilon^M(y_{t})+\int^{s}_{{t}}(\partial_x\Upsilon^M(y_\sigma),F(X^{\gamma_t,u}_\sigma,u(\sigma))_Hd\sigma,
\end{eqnarray*}
where
 $y(s)=X^{\gamma_t,u}(s)-e^{A(s-t)}\eta_t(t),\  t\leq s\leq T$ and $y(s)=\gamma_t(s)-\eta_t(s),\  0\leq s<t$.
  That is (\ref{jias51}).
The proof is now complete. \ \ $\Box$
\begin{remark}\label{remarks}
 \begin{description}
  \rm{
 \item{(i)}      Since $||\cdot||_{0}^2$ is not belongs to $C^{1}(\hat{\Lambda})$, then, for every $a_{\hat{t}}\in \Lambda$, $||\gamma_t-a_{\hat{t},t,A}||_0^2$ cannot  appear  as an auxiliary functional in the proof of the
    uniqueness and stability of viscosity solutions. However, by the above lemma, we can replace $||\gamma_t-a_{\hat{t},t,A}||_{0}^2$
   with its equivalent functional  $\Upsilon^2(\gamma_t-a_{\hat{t},t,A})$. Therefore, we can get the uniqueness result of viscosity solutions without the the assumption $(A.4)$ on page 231 of {\cite{li}}.
 \item{(ii)}
       It follows from (\ref{s0}) that $\overline{\Upsilon}^2(\cdot,\cdot)$ is a gauge-type function.  We can apply it to Lemma \ref{theoremleft} to
get a maximum of a perturbation of the auxiliary functional in the proof of uniqueness.
% For notational simplicity,
%we use $S_1(\cdot)$ to denote $S_1(\cdot,\mathbf{0})$, which will be applied   to prove the uniqueness of viscosity solutions.
  }
\end{description}
\end{remark}
%\begin{remark}\label{remarks}
%   Since $||\gamma_t-a_{\hat{t}}||_0^2$ is not belongs to $C^{1}(\Lambda^{\hat{t}})$, it cannot  appear  as an auxiliary functional in the proof of the
%  uniqueness of viscosity solutions. However, by the above lemma and Lemma \ref{theoremS}, we can replace $||\gamma_t-a_{\hat{t}}||_0^2$   with its equivalent functional  $\Upsilon^2(\gamma_t-a_{\hat{t},t,A})\in C^1(\Lambda^{\hat{t}})$.
%\end{remark}

%%%%%%%%%%%%%%%%%%%%%%%%%%%%%%%%%%%%%%%%%%%%%%%%%%%%%%%%%%%%%%%
\section{Viscosity solutions to  PHJB equations: Existence theorem.}%Construction of random invariant manifold

\par
In this section, we consider the  first order path-dependent  Hamilton-Jacobi-Bellman
                   (PHJB) equation (\ref{hjb1}). As usual, we start with classical solutions.
                                  % In this section, we consider the following second order path-dependent HJB (PHJB) equation:
%\begin{eqnarray}\label{hjb}
%\cases{
%           \partial_tV(\gamma_t)+{\mathbf{H}}(\gamma_t,V(\gamma_t),\partial_xV(\gamma_t),\partial_{xx}V(\gamma_t))= 0,\ \ \  t\in
%                               [0,T),\ \  \gamma_t\in {\Lambda},\cr
%  V(\gamma_T)=\phi(\gamma_T), \ \ \ \gamma_T\in {\Lambda}_T;\cr
%}
%\end{eqnarray}
%      where
%\begin{eqnarray*}
%                                {\mathbf{H}}(\gamma_t,r,p,l)&=&\sup_{u\in{
%                                         {U}}}[
%                        (p,F(\gamma_t,u))_{R^d}+\frac{1}{2}\mbox{tr}[ l G(\gamma_t,u)G^\top(\gamma_t,u)]\\
%                        &&\ \ \ \ \ +q(\gamma_t,r,G^\top(\gamma_t,u)p,u)],  \ \ (\gamma_t,r,p,l)\in {\Lambda}\times R\times R^d\times \Gamma(R^{d}).
%\end{eqnarray*}
%Here we let  $G^\top$ the transpose of the matrix $G$,  $\Gamma(R^d)$  the set of all $(d\times d)$ symmetric matrices and  $(\cdot,\cdot)_{R^d}$ the scalar product of $R^d$.
\par
\begin{definition}\label{definitionccc}     (Classical solution)
              A functional $v\in C^{1}({\Lambda})$  with  $A^*\partial_xv\in C^0(\Lambda)$    is called a classical solution to the PHJB equation (\ref{hjb1}) if it satisfies
              the PHJB¡¡  equation (\ref{hjb1}) point-wisely.
 \end{definition}
% \par
% \begin{definition}\label{definitionusc}     We define
 %\begin{eqnarray*}
%                           &&USC_*(\tilde{\Lambda}^M):=\{u:\tilde{\Lambda}\rightarrow R\  \mbox{such that} \\
%                           &&(i)  \ \mbox{For each fixed} \ \omega_{\hat{t}}\in \tilde{\Lambda}, \bar{u}(t,x):=u((\omega^x_{\hat{t}})_{\hat{t},t})
%                           \ \mbox{is a}\ USC\mbox{-function of } \ (t, x+\omega(\hat{t}))\in [0,T-\hat{t}]\times Q;\\
%                           &&(ii) \ \mbox{For each}\ \omega_t\in \tilde{\Lambda}^M \ \mbox{with}\ t_i\uparrow t,  \limsup_{i\rightarrow\infty}u(\omega_{t_i})
%                           \leq \sup_{|x+\omega(t)|\leq M}u(\omega_{t}^x)\}.
 %\end{eqnarray*}
%              A functional $v\in C^{1,2}(\tilde{\Lambda})$       is called a classical solution to the path-dependent HJB equation (\ref{hjb}) if it satisfies
%              the path-dependent HJB¡¡equation (\ref{hjb}) point-wisely.
% \end{definition}
        \par
        We will prove that the value functional $V$ defined by (\ref{value1}) is a viscosity solution of PHJB equation (\ref{hjb1}).
Define
$$
\Phi=\{\varphi\in C^1(\Lambda)| %\phi \ \mbox{is weakly sequentially lower semicontinuous},
A^*\partial_x\varphi\in C^0(\Lambda)\};
$$
\begin{eqnarray*}
{\mathcal{G}}_t&=&\bigg{\{}g\in C^0(\Lambda^t)| \exists \  h\in C^1([0,T]\times R), \delta_i>0, \gamma^i_{t_i}\in \Lambda, t_i\leq t, i=0,1,2,\ldots, N>0, \ \mbox{with} \ \nabla_xh \geq0,\\
&& \ ||\gamma^i_{t_i}||_0\leq N, \sum_{i=0}^{\infty}\delta_i\leq N \ \mbox{such that}\
                   g(\gamma_s)=h(s,\Upsilon^2(\gamma_s))+\sum_{i=0}^{\infty}\delta_i\overline{\Upsilon}^2(\gamma_s-\gamma^i_{t_i,s,A}), \forall\ \gamma_s\in \Lambda^t
\bigg{\}}.
\end{eqnarray*}
For notational simplicity, if $g\in {\mathcal{G}}_t$,
we use $\partial_tg(\gamma_s)$ and $\partial_xg(\gamma_s)$ to denote $h_t(s,\Upsilon^2(\gamma_s))+2\sum_{i=0}^{\infty}\delta_i(s-t_i)$ and $\nabla_xh(s,\Upsilon^2(\gamma_s))\partial_x\Upsilon^2(\gamma_s)+\sum_{i=0}^{\infty}\delta_i\partial_x\Upsilon^2(\gamma_s-\gamma^i_{t_i,s,A})$, respectively.
%\begin{eqnarray*}
%{\mathcal{G}}_t&=&\{g\in C^1(\Lambda^t)|  (A_{\mu}\gamma_s(s),\partial_xg(\gamma_s)\leq0\ \mbox{for all}\ \gamma_s\in \Lambda^t\ \mbox{and}\ \mu>0\}.
%\exists \  g_i\in C^1(R), \gamma^i_{t_i}\in \Lambda, t_i\leq t, i=0,1,2,\ldots, M>0, \ \mbox{with}\ g'_i\geq0, \\
%&&\ ||\gamma^i_{t_i}||_0\leq M \ \mbox{such that}\
%                   g(\gamma_s)=h(s)\sum_{i=0}^{\infty}\frac{1}{2^i}g_i(\Upsilon^2(\gamma_s-\gamma^i_{t_i,s,A})), \forall\ \gamma_s\in \Lambda^t
%\}.
%\end{eqnarray*}
\par
                     Now we can give the following definition for  viscosity solutions.
\begin{definition}\label{definition4.1} \ \
 $w\in C^0({\Lambda})$ is called a
                             viscosity subsolution (resp., supersolution)
                             to  (\ref{hjb1}) if the terminal condition,  $w(\gamma_T)\leq \phi(\gamma_T)$(resp., $w(\gamma_T)\geq \phi(\gamma_T)$),
                             $\gamma_T\in {\Lambda}_T$ is satisfied, and for any $\varphi\in \Phi$ and $g\in {\cal{G}}_t$ with $t\in [0,T)$, whenever the function $w-\varphi-g$  (resp.,  $w+\varphi+g$) satisfies
$$
                         0=({w}-\varphi-g)(\gamma_t)=\sup_{\eta_s\in \Lambda^t}
                         ({w}- \varphi-g)(\eta_s),
$$
$$
                       (\mbox{resp.,}\ \
                         0=({w}+\varphi+g)(\gamma_t)=\inf_{\eta_s\in \Lambda^t}
                         ({w}+\varphi+g)(\eta_s),)
$$
                   %where $\hat{\gamma}^\mu_{s_\mu}\in {\cal{C}}^\alpha_{\mu,M_0}$, $s_\mu\in [0,T-\Delta]$ and $|\hat{\gamma}^\mu_{s_\mu}(s_\mu)|<
%                   \frac{M_0}{2}$,
              where $\gamma_t\in \Lambda$,      we have
\begin{eqnarray*}
                          \partial_t\varphi(\gamma_t)+\partial_tg(\gamma_t)+(A^*\partial_x\varphi(\gamma_t),\gamma_t(t))_H
                           +{\mathbf{H}}(\gamma_t,\partial_x\varphi(\gamma_t)+\partial_xg(\gamma_t))\geq0,
\end{eqnarray*}
\begin{eqnarray*}
                          (\mbox{resp.,}\ -\partial_t\varphi(\gamma_t)-\partial_tg(\gamma_t)-(A^*\partial_x\varphi(\gamma_t),\gamma_t(t))_H
                           +{\mathbf{H}}(\gamma_t,-\partial_x\varphi(\gamma_t)-\partial_xg(\gamma_t))
                          \leq0.)
\end{eqnarray*}
                                $w\in C^0({\Lambda})$ is said to be a
                             viscosity solution to (\ref{hjb1}) if it is
                             both a viscosity subsolution and a viscosity
                             supersolution.
\end{definition}
%\begin{remark}\label{remarkv}
%  %\begin{description}
%%  \rm{
% % \item{(i)}  A viscosity solution $V$ of the PHJB equation (\ref{hjb})  is a
%%                       classical solution (See Definition \ref{definitionccc}) if it further lies in $C^{1,2}_{p,C_0}({\Lambda})$ for some constants $p,C_0>0$.
% %\item{(i)}    In the classical uniqueness proof of viscosity solution to second order HJB equation in infinite
%%                dimensions, the weak compactness of closed balls in separable
%%                Hilbert spaces is used  (See  \cite{swi}). In our case, the  path-dependent HJB
%%                equation is defined on space  ${\Lambda}$, which
%%                does not  have weak compactness.  In order to obtain  the
%%uniqueness of viscosity solutions, the compact subset  sequence $\{{\cal{C}}^\alpha_{\mu,M_0}\}_{\mu>0}$ are used to introduce the new concept of viscosity  solution.
%%\item{(ii)}
%  Assume that  the coefficients $F(\gamma_t,u)=\overline{F}(t,\gamma_t(t),u),
%                     \ G(\gamma_t,u)=\overline{G}(t,\gamma_t(t),u)$,  $ q(\gamma_t,$ $y,z,u)=\overline{q}(t,\gamma_t(t),y,z,u),
%                     \phi(\gamma_T)=\overline{\phi}(\gamma_T(T)), \
%                     (\gamma,y,z,u) \in {\Lambda}_T\times R\times R^d\times U$.
%           Let the  function $V:[0,T]\times R^d\rightarrow R$ be a
%             viscosity solution to (\ref{hjb1}) as a  functional of  $V(\gamma_t):{\Lambda}\rightarrow R$. Then, $V$ is also a classical viscosity
%solution as a function of the  state.
%%}
%%\end{description}
%\end{remark}
\begin{theorem}\label{theoremvexist} \ \
                          Suppose that Hypothesis \ref{hypstate}  holds. Then the value
                          functional $V$ defined by (\ref{value1}) is a
                          viscosity solution to (\ref{hjb1}).
\end{theorem}

 {\bf  Proof}. \ \         First, Let  $\varphi\in \Phi$ and $g\in {\cal{G}}_t$ with $t\in [0,T)$
 such that
$$
                         0=(V-\varphi-g)(\gamma_t)=\sup_{\eta_s\in \Lambda^t}
                         (V- \varphi-g)(\eta_s),
$$
 where $\gamma_t\in \Lambda$.
 Thus, for fixed $u\in U$, by the DPP (Theorem \ref{theorem3.3}), we obtain that, for all $t\leq t+\delta\leq T$,
 \begin{eqnarray}\label{4.9}
                            (\varphi+g)(\gamma_t)&=&V(\gamma_t)
                           \leq \int^{t+\delta}_{t}
                              q(X_\sigma^{\gamma_t,u},u)d\sigma
                              +V(X_{t+\delta}^{\gamma_t,u})\nonumber\\
                              &
                              \leq& \int^{t+\delta}_{t}
                              q(X_\sigma^{\gamma_t,u},u)d\sigma
                              +(\varphi+g)(X_{t+\delta}^{\gamma_t,u}).
\end{eqnarray}
Applying Theorem  \ref{theoremito} and Lemma \ref{theoremS000}, we show that
  \begin{eqnarray*}
                 0  &\leq&\lim_{\delta\rightarrow0}\bigg{[}\frac{1}{\delta}\int^{t+\delta}_{t}
                              q(X_\sigma^{\gamma_t,u},u)d\sigma
                            +\frac{1}{\delta}{[}{{\varphi}}(X_{t+\delta}^{\gamma_t,u})-{{\varphi}} (\gamma_t)]+\frac{1}{\delta}{[}g(X_{t+\delta}^{\gamma_t,u})-g(\gamma_t){]}\bigg{]}\\
                            &\leq&q(\gamma_t,u)+ \partial_t\varphi(\gamma_t)+ \partial_tg(\gamma_t)+(A^*\partial_x\varphi(\gamma_t),\gamma_t(t))_H+(\partial_{x}\varphi(\gamma_t)+\partial_{x}g(\gamma_t),
                     F(\gamma_t,u))_H.
\end{eqnarray*}
   Taking the minimum in $u\in U$, %and letting   $\mu\rightarrow+\infty$,
   we have that $V$ is a viscosity subsolution of (\ref{hjb1}).
%$$
%0
%    \leq\partial_t{\varphi}(\hat{\gamma}_{t})+(A^*\partial_x\varphi(\gamma_t),\gamma_t(t))_H
%                           +{\mathbf{H}}(\hat{\gamma}_{t},\partial_x{\varphi}(\hat{\gamma}_{t})).
%$$
\par
      Next,   Let  $\varphi\in \Phi$ and $g\in {\cal{G}}_t$ with $t\in [0,T)$
 such that
$$
                         0=(V+\varphi+g)(\gamma_t)=\inf_{\eta_s\in \Lambda^t}
                         (V+\varphi+g)(\eta_s),
$$
 where $\gamma_t\in \Lambda$.
 Then, for any $\varepsilon>0$,  by  the DPP (Theorem \ref{theorem3.3}), one can find a control  ${u}^\varepsilon(\cdot)\equiv u^{\varepsilon,\delta}(\cdot)\in {\cal{U}}[t,T]$ such
   that, for any $t\leq t+\delta\leq T$,
\begin{eqnarray*}
    \varepsilon\delta
    &\geq& \int^{t+\delta}_{t}
                              q(X_\sigma^{\gamma_t,{u}^\varepsilon},{u}^\varepsilon(\sigma))d\sigma
                              +V(X_{t+\delta}^{\gamma_t,{u}^\varepsilon})-V(\gamma_t)\\
                              &
                              \geq& \int^{t+\delta}_{t}
                              q(X_\sigma^{\gamma_t,{u}^\varepsilon},{u}^\varepsilon(\sigma))d\sigma
                              -(\varphi+g)(X_{t+\delta}^{\gamma_t,{u}^\varepsilon})+(\varphi+g)(\gamma_t).
\end{eqnarray*}
                 Then, by Theorem  \ref{theoremito} and Lemma \ref{theoremS000},   we obtain that
\begin{eqnarray*}
                           \varepsilon&\geq& \frac{1}{\delta}\int^{t+\delta}_{t}
                              q(X_\sigma^{\gamma_t,{u}^\varepsilon},{u}^\varepsilon(\sigma))d\sigma
                                -\frac{\varphi(X_{t+\delta}^{\gamma_t,{u}^\varepsilon})-\varphi(\gamma_t)}{\delta}-\frac{g(X_{t+\delta}^{\gamma_t,{u}^\varepsilon})-g(\gamma_t)}{\delta}
                                 \\
           &\geq& -\partial_t\varphi(\gamma_t)-\partial_tg(\gamma_t)-(A^*\partial_x(\gamma_t),\gamma_t(t))_H+\frac{1}{\delta}\int^{t+\delta}_{t}q(\gamma_t,{u}^\varepsilon(\sigma)) \\
           &&~~~~~~-(\partial_{x}\varphi(\gamma_t)+\partial_{x}g(\gamma_t),
                     F(\gamma_t,{u}^\varepsilon(\sigma)))_H d\sigma+o(1)\\
            &\geq& -\partial_t\varphi(\gamma_t)-\partial_tg(\gamma_t)-(A^*\partial_x(\gamma_t),\gamma_t(t))_H
            +\inf_{u\in U}[q(\gamma_t,{u}) -(\partial_{x}\varphi(\gamma_t)+\partial_{x}g(\gamma_t),
                     F(\gamma_t,{u}))_H]+o(1).
\end{eqnarray*}
Letting $\delta\downarrow 0$ and  $\varepsilon\rightarrow0$, we show that
$$
            0\geq -\partial_t\varphi(\gamma_t)-\partial_tg(\gamma_t)-(A^*\partial_x(\gamma_t),\gamma_t(t))_H
            +H(\gamma_t,-\varphi(\gamma_t)-g(\gamma_t)).
$$
    Therefore, $V$ is also a  viscosity supsolution of equation (\ref{hjb1}).  This completes the proof.\ \ $\Box$
 \par
 \par
Now, let us give   the result of  classical  solutions, which show the consistency of viscosity solutions.
 \begin{theorem}\label{theorem3.2}
                      Let $V$ denote the value functional  defined by (\ref{value1}). If  $V\in C^{1}({\Lambda})$ and $A^*\partial_xV\in C^0(\Lambda)$, then
                 $V$ is a classical solution of  equation (\ref{hjb1}).
\end{theorem}
{\bf  Proof }. \ \
First, using the  definition of $V$ yields  $V(\gamma_T)=\phi(\gamma_T)$ for all $\gamma_T\in \Lambda_T$. Next, for   fixed $(t,\gamma_t,u)\in [0,T)\times{\Lambda}\times U$,
                  from  the DPP (Theorem \ref{theorem3.3}), we obtain the following result:
 \begin{eqnarray}\label{4.9000}
                           && 0\leq \int_{t}^{t+\delta}q(X^{\gamma_t,u}_\sigma,u)d\sigma
                                +V(X^{\gamma_t,u}_{t+\delta})
                           -V(\gamma_t),\ \ 0<\delta< T-t.
\end{eqnarray}
%Thus
%\begin{eqnarray*}
%                 0&\geq& E\bigg{[}\int_{t}^{t+\delta}q(X^{{\gamma}_{{t}},{u}}_s,Y^{{\gamma}_{{t}},{u}}(s),Z^{{\gamma}_{{t}},
%               {u}}(s),{u})ds
%                                +V(X^{{\gamma}_{{t}},u}_{{{t}}+\delta})\bigg{]}
%                 -V(\gamma_t).
%\end{eqnarray*}
               By Theorem  \ref{theoremito},   the  inequality above implies that
               \begin{eqnarray*}
    0&\leq&\lim_{\delta\rightarrow 0^+}\frac{1}{\delta}\bigg{[}\int_{t}^{t+\delta}q(X^{\gamma_t,u}_\sigma,u)d\sigma+V(X^{\gamma_t,u}_{t+\delta})-V(\gamma_t)\bigg{]}\\
                       &=&\partial_t V(\gamma_t)+(A^*\partial_xV(\gamma_t),\gamma_t(t))_H +(
                        F(\gamma_t,u),\partial_xV(\gamma_t))_{H}+q(\gamma_t,u).
\end{eqnarray*}
     Taking the minimum in $u\in U$, we have that
\begin{eqnarray}\label{3.14}
    0\leq\partial_t V(\gamma_t)+(A^*\partial_xV(\gamma_t),\gamma_t(t))_H + H(\gamma_t,\partial_x V(\gamma_t)).
\end{eqnarray}
          On the other hand,  let $(t,\gamma_t)\in [0,T)\times\Lambda$ be fixed. Then, by (\ref{3.7}) and Theorem  \ref{theoremito},  there exists an
                         $\tilde{u}\equiv u^{\varepsilon,\delta}\in {\cal{U}}[t,T]$ for any $\varepsilon>0$ and $0<\delta<T-t$ such that
%\begin{eqnarray*}
%    &&-\varepsilon\delta \leq G^{\gamma_t,{\tilde{u}}}_{t,t+\delta}[V(X^{{\gamma}_{{t}},{\tilde{u}}}_{{{t}}+\delta})]
%                           -V(\gamma_t).
%\end{eqnarray*}
%Thus we have
\begin{eqnarray*}
                \varepsilon\delta &\geq& \int_{t}^{t+\delta}q(X^{{\gamma}_{{t}},{\tilde{u}}}_s,{\tilde{u}}(s))ds
                                +V(X^{{\gamma}_{{t}},{\tilde{u}}}_{{{t}}+\delta})
                 -V(\gamma_t)\\
                  &=& {\partial_t}V(\gamma_t)\delta+(A^*\partial_xV(\gamma_t),\gamma_t(t))_H\delta +\int_{t}^{t+\delta}q(\gamma_t,\tilde{u}(\sigma))+({\partial_{x}V(\gamma_t)},
                        F(\gamma_t,{\tilde{u}}(\sigma))_Hd\sigma+o(\delta)\\
                     &\geq& {\partial_t}V(\gamma_t)\delta+(A^*\partial_xV(\gamma_t),\gamma_t(t))_H\delta+H(\gamma_t,\partial_x V(\gamma_t))\delta+o(\delta).
\end{eqnarray*}
Then, dividing through by $\delta$ and letting $\delta\rightarrow0^+$, we
obtain that
$$
                       \varepsilon\geq \partial_tV(\gamma_t)+(A^*\partial_xV(\gamma_t),\gamma_t(t))_H+{\mathbf{H}}(\gamma_t,\partial_{x}V(\gamma_t)).
$$
               The desired result  is obtained by combining the inequality given above with (\ref{3.14}). \ \ $\Box$
                                \par
We conclude this section with   the stability of viscosity solutions.
\begin{theorem}\label{theoremstability}
                      Let $F,q,\phi$ satisfy  Hypothesis \ref{hypstate}, and $v\in C^0(\Lambda)$. Assume
                       \item{(i)}      for any $\varepsilon>0$, there exist $F^\varepsilon,  q^\varepsilon, \phi^\varepsilon$ and $v^\varepsilon\in C^0(\Lambda)$ such that  $F^\varepsilon, q^\varepsilon, \phi^\varepsilon$ satisfy  Hypothesis \ref{hypstate} and $v^\varepsilon$ is a viscosity subsolution (resp., supsolution) of equation (\ref{hjb1}) with generators $F^\varepsilon,  q^\varepsilon, \phi^\varepsilon$;
                           \item{(ii)} as $\varepsilon\rightarrow0$, $(F^\varepsilon,  q^\varepsilon, \phi^\varepsilon,v^\varepsilon)$ converge to
                           $(F,  q, \phi, v)$ uniformly in the following sense: %for any $(t,\gamma_t)\in [0,T]\times \Lambda$, there exists $\delta$ such that
\begin{eqnarray}\label{sss}
                         \lim_{\varepsilon\rightarrow0}\sup_{(t,{\gamma_t},u)\in [0,T]\times\Lambda\times U}\sup_{\eta_T\in \Lambda_T}[(|F^\varepsilon-F|
                         +|q^\varepsilon-q|)({\gamma}_t,u)+|\phi^\varepsilon-\phi|(\eta_T)+|v^\varepsilon-v|({\gamma}_t)]=0.
\end{eqnarray}
                   Then $v$ is a viscosity subsoluiton (resp., supersolution) of equation (\ref{hjb1}) with generators $F,q,\phi$.
\end{theorem}
{\bf  Proof }. \ \ Without loss of generality, we shall only prove the viscosity subsolution property.
First,  from $v^{\varepsilon}$ is a viscosity subsolution of  equation (\ref{hjb1}) with generators $F^{\varepsilon},   q^{\varepsilon}, \phi^{\varepsilon}$, it follows that
 $$
                            v^{\varepsilon}(\gamma_T)\leq \phi^{\varepsilon}(\gamma_T),\ \ \gamma_T\in \Lambda_T.
$$
Letting $\varepsilon\rightarrow0$, we  have
 $$
                            v(\gamma_T)\leq \phi(\gamma_T),\ \ \gamma_T\in \Lambda_T.
$$
Next,     Let  $\varphi\in \Phi$ and $g\in {\cal{G}}_{\hat{t}}$ with $\hat{t}\in [0,T)$
 such that
$$
                         0=(V-\varphi-g)(\hat{\gamma}_{\hat{t}})=\sup_{\eta_s\in \Lambda^{\hat{t}}}
                         (V- \varphi-g)(\eta_s),
$$
 where $\hat{\gamma}_{\hat{t}}\in \Lambda$. % we let   $\varphi\in J^+(\hat{\gamma}_{\hat{t}}, v)$ with
%  $(\hat{t},\hat{\gamma}_{\hat{t}})\in [0,T)\times\Lambda$.
 Denote $g_{1}(\gamma_t):=g(\gamma_t)+\overline{\Upsilon}^2(\gamma_t,\hat{\gamma}_{{\hat{t}}})$
 %+|\gamma_t(t)-\hat{\gamma}_{{\hat{t}}}(\hat{t})|^2$
 for all
 $(t,\gamma_t)\in [\hat{t},T]\times\Lambda$. %By Lemma \ref{theoremS},
 Then we have  $g_{1}\in {\cal{G}}_{\hat{t}}$.  Define a sequence of positive numbers $\{\delta_i\}_{i\geq0}$  by %$\delta_0=\beta$ and
        $\delta_i=\frac{1}{2^i}$ for all $i\geq0$.
 For every $\varepsilon>0$, since $v^{\varepsilon}-\varphi-g_1$ is a  upper semicontinuous functional  and $\overline{\Upsilon}^2(\cdot,\cdot)$ is a gauge-type function, from Lemma \ref{theoremleft} it follows that,
  for every  $(t_0,\gamma^0_{t_0})\in [\hat{t},T]\times \Lambda^{\hat{t}}$ satisfy
 $$
(v^{\varepsilon}-\varphi-g_1)(\gamma^0_{t_0})\geq \sup_{(s,\gamma_s)\in [\hat{t},T]\times \Lambda^{\hat{t}}}(v^{\varepsilon}-\varphi-g_1)(\gamma_s)-\varepsilon,\
\    \mbox{and} \ \ (v^{\varepsilon}-\varphi-g_1)(\gamma^0_{t_0})\geq (v^{\varepsilon}-\varphi-g_1)(\hat{\gamma}_{\hat{t}}),
 $$
  there exist $(t_{\varepsilon},{\gamma}^{\varepsilon}_{t_{\varepsilon}})\in [\hat{t},T]\times \Lambda^{\hat{t}}$ and a sequence $\{(t_i,\gamma^i_{t_i})\}_{i\geq1}\subset [\hat{t},T]\times \Lambda^{\hat{t}}$ such that
  \begin{description}
        \item{(i)} $\overline{\Upsilon}^2(\gamma^0_{t_0},{\gamma}^{\varepsilon}_{t_{\varepsilon}})\leq {\varepsilon}$,  $\overline{\Upsilon}^2(\gamma^i_{t_i},{\gamma}^{\varepsilon}_{t_{\varepsilon}})\leq \frac{\varepsilon}{2^i}$ and $t_i\uparrow t_{\varepsilon}$ as $i\rightarrow\infty$,
        \item{(ii)}  $(v^{\varepsilon}-\varphi-g_1)({\gamma}^{\varepsilon}_{t_{\varepsilon}})-\sum_{i=0}^{\infty}\frac{1}{2^i}\overline{\Upsilon}^2(\gamma^i_{t_i},{\gamma}^{\varepsilon}_{t_{\varepsilon}})\geq (v^{\varepsilon}-\varphi-g_1)(\gamma^0_{t_0})$, and
        \item{(iii)}  $(v^{\varepsilon}-\varphi-g_1)(\gamma_s)-\sum_{i=0}^{\infty}\frac{1}{2^i}\overline{\Upsilon}^2(\gamma^i_{t_i},\gamma_s)
            <(v^{\varepsilon}-\varphi-g_1)({\gamma}^{\varepsilon}_{t_{\varepsilon}})-\sum_{i=0}^{\infty}\frac{1}{2^i}\overline{\Upsilon}^2(\gamma^i_{t_i},{\gamma}^{\varepsilon}_{t_{\varepsilon}})$ for all $(s,\gamma_s)\in [t_{\varepsilon},T]\times \Lambda^{t_{\varepsilon}}\setminus \{(t_{\varepsilon},{\gamma}^{\varepsilon}_{t_{\varepsilon}})\}$.

        \end{description}
  %
%
%
%
%  there exists $({{t}_{\varepsilon}},{\gamma}^{\varepsilon}_{{t}_{\varepsilon}})\in [\hat{t},T]\times {\Lambda}_{M_0}(\hat{\gamma}_{\hat{t}})$ such that
%$$
%                    (v^{\varepsilon}- {{\varphi_{1}}})( \hat{\gamma}_{\hat{t}})\leq   (v^{\varepsilon}- {{\varphi_{1}}})({\gamma}^{\varepsilon}_{{t}_{\varepsilon}})=\sup_{\gamma_s\in {\Lambda}_{M_0}({\gamma}^{\varepsilon}_{{t}_{\varepsilon}})}
%                         (v^{\varepsilon}- {{\varphi_{1}}})(
%                         \gamma_s).
%$$
We claim that
\begin{eqnarray}\label{gamma}
d_\infty({\gamma}^{\varepsilon}_{{t}_{\varepsilon}},\hat{\gamma}_{\hat{t}})\rightarrow0  \ \ \mbox{as} \ \ \varepsilon\rightarrow0.
\end{eqnarray}
 Indeed, if not,  by (\ref{s0}), we can assume that there exists an $\nu_0>0$
 such
                    that
$$
                 %   |t_\varepsilon-\hat{t}|^2
%  +
  \overline{\Upsilon}^2({\gamma}^{\varepsilon}_{{t}_{\varepsilon}},\hat{\gamma}_{{\hat{t}}})%+|{\gamma}^{\varepsilon}_{{t}_{\varepsilon}}({t}_{\varepsilon})-\hat{\gamma}_{{\hat{t}}}(\hat{t})|_0^8
  \geq\nu_0.
$$
 Thus,  we obtain that
%Now we have  $({\bar{t}},\bar{\gamma}_{\bar{t}})=(\hat{t},\hat{\gamma}_{\hat{t}})$, for if not,
\begin{eqnarray*}
   &&0=(v- {{\varphi}}-g)(\hat{\gamma}_{\hat{t}})= \lim_{\varepsilon\rightarrow0}(v^\varepsilon-\varphi-g_1)(\hat{\gamma}_{\hat{t}})
   \leq \overline{\lim_{\varepsilon\rightarrow0}}\bigg{[}(v^{\varepsilon}-\varphi-g_1)({\gamma}^{\varepsilon}_{t_{\varepsilon}})-\sum_{i=0}^{\infty}\frac{1}{2^i}\overline{\Upsilon}^2(\gamma^i_{t_i},{\gamma}^{\varepsilon}_{t_{\varepsilon}})\bigg{]}\\
   %&=&\overline{\lim_{\varepsilon\rightarrow0}}\bigg{[}(v^\varepsilon-{{\varphi}}-g)({\gamma}^{\varepsilon}_{{t}_{\varepsilon}})%-|t_\varepsilon-\hat{t}|^2
%  -\overline{\Upsilon}^2({\gamma}^{\varepsilon}_{{t}_{\varepsilon}},\hat{\gamma}_{{\hat{t}}})-%|{\gamma}^{\varepsilon}_{{t}_{\varepsilon}}({t}_{\varepsilon})-\hat{\gamma}_{{\hat{t}}}(\hat{t})|_8^8
%  \sum_{i=0}^{\infty}\frac{1}{2^i}\overline{\Upsilon}^2(\gamma^i_{t_i},{\gamma}^{\varepsilon}_{t_{\varepsilon}})\bigg{]}\\
   &\leq&\overline{\lim_{\varepsilon\rightarrow0}}{[}(v-{{\varphi}}-g)({\gamma}^{\varepsilon}_{{t}_{\varepsilon}})+(v^\varepsilon-v)({\gamma}^{\varepsilon}_{{t}_{\varepsilon}})
     %-\sum_{i=0}^{\infty}\frac{1}{2^i}\overline{\Upsilon}^2(\gamma^i_{t_i},{\gamma}^{\varepsilon}_{t_{\varepsilon}})
     {]}-\nu_0\leq (v- {{\varphi}}-g)(\hat{\gamma}_{\hat{t}})-\nu_0=-\nu_0,
\end{eqnarray*}
 contradicting $\nu_0>0$.  We notice that, by  the property (i) of $(t_{\varepsilon},{\gamma}^{\varepsilon}_{t_{\varepsilon}})$,
  \begin{eqnarray*}
  2\sum_{i=0}^{\infty}\frac{1}{2^i}({t_{\varepsilon}}-{t}_{i})
  \leq2\sum_{i=0}^{\infty}\frac{1}{2^i}\bigg{(}\frac{\varepsilon}{2^i}\bigg{)}^{\frac{1}{2}}\leq 4\varepsilon^{\frac{1}{2}};
    \end{eqnarray*}
    \begin{eqnarray*}
    |\partial_x{\Upsilon}^2({\gamma}^{\varepsilon}_{{t}_{\varepsilon}}-\hat{\gamma}_{{\hat{t}},{t}_{\varepsilon},A})|\leq 4|e^{({t}_{\varepsilon}-\hat{t})A}\hat{\gamma}_{{\hat{t}}}(\hat{t})-{\gamma}^{\varepsilon}_{{t}_{\varepsilon}}({t}_{\varepsilon})|;
    \end{eqnarray*}
    and
     \begin{eqnarray*}
  \bigg{|}\partial_x\sum_{i=0}^{\infty}\frac{1}{2^i}
                      \Upsilon^2({\gamma}^{\varepsilon}_{t_{\varepsilon}}-\gamma^i_{t_i,t_{\varepsilon},A})
                      \bigg{|}
                      &\leq&4\sum_{i=0}^{\infty}\frac{1}{2^i}|e^{(t_{\varepsilon}-t_i)A}\gamma^i_{t_i}({t}_{i})-{\gamma}^{\varepsilon}_{t_{\varepsilon}}(t_{\varepsilon})|
                     \leq4\sum_{i=0}^{\infty}\frac{1}{2^i}\bigg{(}\frac{\varepsilon}{2^i}\bigg{)}^{\frac{1}{2}}
                      \leq8{\varepsilon}^{\frac{1}{2}}.
                        \end{eqnarray*}
 Then for any $\varrho>0$, by (\ref{sss}) and (\ref{gamma}), there exists $\varepsilon>0$ small enough such that
$$
            \hat{t}\leq {t}_{\varepsilon}< T,  \        \  %|{\gamma}^{\varepsilon}_{{t}_{\varepsilon}}({t}_{\varepsilon})|< \frac{M_0}{2},\ \
             2|{t}_{\varepsilon}-\hat{t}|+2\sum_{i=0}^{\infty}\frac{1}{2^i}({t_{\varepsilon}}-{t}_{i})+|\partial_t{\varphi}({\gamma}^{\varepsilon}_{{t}_{\varepsilon}})-\partial_t{\varphi}(\hat{\gamma}_{\hat{t}})|
             +|\partial_tg({\gamma}^{\varepsilon}_{{t}_{\varepsilon}})-\partial_tg(\hat{\gamma}_{\hat{t}})|\leq \frac{\varrho}{4},$$
$$
|(A^*\partial_x\varphi({\gamma}^{\varepsilon}_{{t}_{\varepsilon}}),{\gamma}^{\varepsilon}_{{t}_{\varepsilon}}({t}_{\varepsilon}))_H-(A^*\partial_x{\varphi}(\hat{\gamma}_{\hat{t}}),\hat{\gamma}_{\hat{t}}(\hat{t}))_H| \leq \frac{\varrho}{4}, \ \mbox{and}\  |I|+|II|\leq \frac{\varrho}{4},
$$
where
\begin{eqnarray*}
I&=&{\mathbf{H}}^{\varepsilon}({\gamma}^{\varepsilon}_{{t}_{\varepsilon}},
                           \partial_x\varphi({\gamma}^{\varepsilon}_{{t}_{\varepsilon}})+ \partial_xg_2({\gamma}^{\varepsilon}_{{t}_{\varepsilon}}))
                           -{\mathbf{H}}({\gamma}^{\varepsilon}_{{t}_{\varepsilon}},
                           \partial_x\varphi({\gamma}^{\varepsilon}_{{t}_{\varepsilon}})+\partial_xg_2({\gamma}^{\varepsilon}_{{t}_{\varepsilon}})),\\
II&=&{\mathbf{H}}({\gamma}^{\varepsilon}_{{t}_{\varepsilon}},
                           \partial_x\varphi({\gamma}^{\varepsilon}_{{t}_{\varepsilon}})+\partial_xg_2({\gamma}^{\varepsilon}_{{t}_{\varepsilon}}))
                       -{\mathbf{H}}(\hat{\gamma}_{\hat{t}},\partial_x{\varphi}(\hat{\gamma}_{\hat{t}})+\partial_xg(\hat{\gamma}_{\hat{t}})),\\
g_2({\gamma}^{\varepsilon}_{{t}_{\varepsilon}})&=&g({\gamma}^{\varepsilon}_{{t}_{\varepsilon}})+\overline{\Upsilon}^2({\gamma}^{\varepsilon}_{{t}_{\varepsilon}}-\hat{\gamma}_{{\hat{t}},{t}_{\varepsilon},A})
+\sum_{i=0}^{\infty}\frac{1}{2^i}\overline{\Upsilon}^2({\gamma}^{\varepsilon}_{t_{\varepsilon}}-\gamma^i_{t_i,t_{\varepsilon},A}),
\\
\partial_xg_2({\gamma}^{\varepsilon}_{{t}_{\varepsilon}})&=&\partial_xg({\gamma}^{\varepsilon}_{{t}_{\varepsilon}})+\partial_x{\Upsilon}^2({\gamma}^{\varepsilon}_{{t}_{\varepsilon}}-\hat{\gamma}_{{\hat{t}},{t}_{\varepsilon},A})
+\sum_{i=0}^{\infty}\frac{1}{2^i}\partial_x{\Upsilon}^2({\gamma}^{\varepsilon}_{t_{\varepsilon}}-\gamma^i_{t_i,t_{\varepsilon},A}),\\
\partial_tg_2({\gamma}^{\varepsilon}_{{t}_{\varepsilon}})&=&\partial_tg({\gamma}^{\varepsilon}_{{t}_{\varepsilon}})+2({t}_{\varepsilon}-\hat{t})
+2\sum_{i=0}^{\infty}\frac{1}{2^i}(t_{\varepsilon}-t_i),
\end{eqnarray*}
and
\begin{eqnarray*}
                                {\mathbf{H}}^{\varepsilon}(\gamma_t,p)=\sup_{u\in{
                                         {U}}}[
                        (p,F^{\varepsilon}(\gamma_t,u))_{H}+q^{\varepsilon}(\gamma_t,u)],  \ \ (t,\gamma_t,p)\in [0,T]\times {\Lambda}\times H.
\end{eqnarray*}
%$$
%                                {\mathbf{H}}^{\varepsilon}(\gamma_t,p)=\inf_{u\in{
%                                         {U}}}[
%                        (p,F^{\varepsilon}(\gamma_t,u))_{R^d}
%                        +q^{\varepsilon}(\gamma_t,u)],  \ \ (\gamma_t,p)\in {\Lambda}\times R^d.
%$$
 Since $v^{\varepsilon}$ is a viscosity subsolution of PHJB equation (\ref{hjb1}) with generators $F^{\varepsilon},   q^{\varepsilon}, \phi^{\varepsilon}$, we have
$$
                          \partial_t\varphi({\gamma}^{\varepsilon}_{{t}_{\varepsilon}})+\partial_tg_2({\gamma}^{\varepsilon}_{{t}_{\varepsilon}})
                           +(A^*\partial_x\varphi({\gamma}^{\varepsilon}_{{t}_{\varepsilon}}),{\gamma}^{\varepsilon}_{{t}_{\varepsilon}}({t}_{\varepsilon}))_H+{\mathbf{H}}^{\varepsilon}({\gamma}^{\varepsilon}_{{t}_{\varepsilon}},
                           \partial_x\varphi({\gamma}^{\varepsilon}_{{t}_{\varepsilon}})+\partial_xg_2({\gamma}^{\varepsilon}_{{t}_{\varepsilon}}))\geq0.
$$
Thus
\begin{eqnarray*}
                       0&\leq&  \partial_t{\varphi}({\gamma}^{\varepsilon}_{{t}_{\varepsilon}})+\partial_tg({\gamma}^{\varepsilon}_{{t}_{\varepsilon}})
                       +2({t}_{\varepsilon}-\hat{t})+2\sum_{i=0}^{\infty}\frac{1}{2^i}({t_{\varepsilon}}-{t}_{i})
                       +(A^*\partial_x\varphi({\gamma}^{\varepsilon}_{{t}_{\varepsilon}}),{\gamma}^{\varepsilon}_{{t}_{\varepsilon}}({t}_{\varepsilon}))_H\\
                       &&+{\mathbf{H}}(\hat{\gamma}_{\hat{t}},\partial_x{\varphi}(\hat{\gamma}_{\hat{t}})+\partial_xg(\hat{\gamma}_{\hat{t}}))+I+II\\
                       &\leq&\partial_t{\varphi}(\hat{\gamma}_{\hat{t}})+\partial_tg(\hat{\gamma}_{\hat{t}})+(A^*\partial_x{\varphi}(\hat{\gamma}_{\hat{t}}),\hat{\gamma}_{\hat{t}}(\hat{t}))_H
                       +{\mathbf{H}}(\hat{\gamma}_{\hat{t}},\partial_x{\varphi}(\hat{\gamma}_{\hat{t}})+\partial_xg(\hat{\gamma}_{\hat{t}}))+\varrho.
\end{eqnarray*}
Letting $\varrho\downarrow 0$, we show that
$$
\partial_t{\varphi}(\hat{\gamma}_{\hat{t}})+\partial_tg(\hat{\gamma}_{\hat{t}})+(A^*\partial_x{\varphi}(\hat{\gamma}_{\hat{t}}),\hat{\gamma}_{\hat{t}}(\hat{t}))_H
+{\mathbf{H}}(\hat{\gamma}_{\hat{t}},\partial_x{\varphi}(\hat{\gamma}_{\hat{t}})+\partial_xg(\hat{\gamma}_{\hat{t}}))\geq0.
$$
Since ${\varphi}\in \Phi$ and $g\in {\cal{G}}_t$ with $t\in [0,T)$  are arbitrary, we see that $v$ is a viscosity subsolution of PHJB equation (\ref{hjb1}) with generators $F,q,\phi$, and thus completes the proof.
\ \ $\Box$

\section{Viscosity solutions to  PHJB equations: Uniqueness theorem.}
\par
             This section is devoted to a  proof of uniqueness of  viscosity
                   solutions to (\ref{hjb1}). This result, together with
                  the results from  the previous section, will be used to characterize
                   the value functional defined by (\ref{value1}).
                   \par
We  now state the main result of this section.
\begin{theorem}\label{theoremhjbm}  Suppose Hypothesis \ref{hypstate}   holds.
                         Let $W_1\in C^0({\Lambda})$ $(\mbox{resp}., W_2\in C^0({\Lambda}))$ be  a viscosity subsolution (resp., supsolution) to equation (\ref{hjb1}) and  let  there exist  constant $L>0$
                          %and $m>0$, %and a  modulus of continuity $\hat{\omega}$,
                        such that, for any  $0\leq t\leq  s\leq T$ and
                        $\gamma_t, \eta_t \in{\Lambda}$,
\begin{eqnarray}\label{w}
                                   |W_1(\gamma_t)|\vee |W_2(\gamma_t)|\leq L (1+||\gamma_t||_0);
                                   \end{eqnarray}
\begin{eqnarray}\label{w1}
                               |W_1(\gamma_{t,s,A})-W_1(\eta_t)|\vee|W_2(\gamma_{t,s,A})-W_2(\eta_t)|\leq
                        L(1+||\gamma_t||_0+||\eta_t||_0)|s-t|+L||\gamma_{t}-\eta_t||_0.
\end{eqnarray}
                   Then  $W_1\leq W_2$.
\end{theorem}

\par
                      Theorems    \ref{theoremvexist} and \ref{theoremhjbm} lead to the result (given below) that the viscosity solution to   PHJB equation given in (\ref{hjb1})
                      corresponds to the value functional  $V$ of our optimal control problem given in (\ref{state1}) and (\ref{value1}).
\begin{theorem}\label{theorem52}\ \
                 Let Hypothesis \ref{hypstate}  hold. Then the value
                          functional $V$ defined by (\ref{value1}) is the unique viscosity
                          solution to (\ref{hjb1}) in the class of functionals satisfying (\ref{w}) and (\ref{w1}).
\end{theorem}
\par
   {\bf  Proof  }. \ \    Theorem \ref{theoremvexist} shows that $V$ is a viscosity solution to equation (\ref{hjb1}).  Thus, our conclusion follows from %Lemma \ref{lemmavaluev},
   Lemma  \ref{lemmaexist111} and Theorem
    \ref{theoremhjbm}.  \ \ $\Box$
%\par
%{\bf  Proof}. \ \
%               By Theorem 4.4, we know that $V$ is a viscosity solution of (3.8). Thus, our conclusion follows from Theorem 3.2 and Theorem 5.1.\ \ $\Box$

\par
  Next, we prove Theorem \ref{theoremhjbm}.   Let $W_1$ be a viscosity solution of PHJB equation (\ref{hjb1}). %We note that it is sufficient to show $W\leq V$ because the inverse of which can be proved in a similar way.
 We  note that for $\delta>0$, the functional
                    defined by $\tilde{W}:=W_1-\frac{\delta}{t}$ is a subsolution
                   for
 \begin{eqnarray*}
\begin{cases}
{\partial_t} \tilde{W}(\gamma_t)+{\mathbf{H}}(\gamma_t, \partial_x \tilde{W}(\gamma_t))
          = \frac{\delta}{t^2}, \ \  \gamma_t\in {\Lambda}, \\
\tilde{W}(\gamma_T)=\phi(\gamma_T).
\end{cases}
\end{eqnarray*}
                As $W_1\leq W_2$ follows from $\tilde{W}\leq W_2$ in
                the limit $\delta\downarrow0$, it suffices to prove
                $W_1\leq W_2$ under the additional assumption given below:
$$
{\partial_t} {W_1}(\gamma_t)+{\mathbf{H}}(\gamma_t, \partial_x {W_1}(\gamma_t))
          \geq c,\ \ c:=\frac{\delta}{T^2}, \ \  \gamma_t\in {\Lambda}.
$$
\par
   {\bf  Proof of Theorem \ref{theoremhjbm} } \ \   The proof of this theorem  is rather long. Thus, we split it into several
        steps.
\par
            $Step\  1.$ Definitions of auxiliary functionals.
            \par
 We only need to prove that $W_1(\gamma_t)\leq W_2(\gamma_t)$ for all $(t,\gamma_t)\in
[T-\bar{a},T)\times
       {\Lambda}$.
        Here,
        $$\bar{a}=\frac{1}{16L}\wedge{T}.$$
         Then, we can  repeat the same procedure for the case
        $[T-i\bar{a},T-(i-1)\bar{a})$.  Thus, we assume the converse result that $(\tilde{t},\tilde{\gamma}_{\tilde{t}})\in [T-\bar{a},T)\times
      {\Lambda}$ exists  such that
        $\tilde{m}:=W_1(\tilde{\gamma}_{\tilde{t}})-W_2(\tilde{\gamma}_{\tilde{t}})>0$.  %Because $\cup_{\mu>0,M_0>0}{\cal{C}}^\alpha_{\mu,M_0}$  is dense in $(\Lambda,d_\infty)$,
%         by (\ref{w}) there exist $\tilde{t}\in [T-\bar{a},T)$ and $\tilde{\gamma}_{\tilde{t}}\in \cup_{\mu>0,M_0>0}{\cal{C}}^\alpha_{\mu,M_0}$
%         such that
%        $W(\tilde{\gamma}_{\tilde{t}})-V(\tilde{\gamma}_{\tilde{t}})>\tilde{m}$.
\par
         Consider that  $\varepsilon >0$ is  a small number such that
 $$
 W_1(\tilde{\gamma}_{\tilde{t}})-W_2(\tilde{\gamma}_{\tilde{t}})-2\varepsilon \frac{\nu T-\tilde{t}}{\nu
 T}\Upsilon^2(\tilde{\gamma}_{\tilde{t}})%-\frac{\varepsilon}{\tilde{t}-T+\bar{a}}
 >\frac{\tilde{m}}{2},
 $$
      and
\begin{eqnarray}\label{5.3}
                          \frac{\varepsilon}{\nu T}\leq\frac{c}{2},
\end{eqnarray}
             where
$$
            \nu=1+\frac{1}{16TL}.
$$
 Next,  we define for any  $(t,\gamma_t,\eta_t)\in (T-\bar{a},T]\times{\Lambda}\times{\Lambda}$,
\begin{eqnarray*}
                 \Psi(\gamma_t,\eta_t)=W_1(\gamma_t)-W_2(\eta_t)-{\beta}\Upsilon^2(\gamma_{t},\eta_{t})-\varepsilon\frac{\nu T-t}{\nu
                 T}(\Upsilon^2(\gamma_t)+\Upsilon^2(\eta_t)).
\end{eqnarray*}
        Define a sequence of positive numbers $\{\delta_i\}_{i\geq0}$  by %$\delta_0=\beta$ and
        $\delta_i=\frac{1}{2^i}$ for all $i\geq0$.
           Since $\Psi$ is a  upper semicontinuous function bounded from above and $\overline{\Upsilon}^2(\cdot,\cdot)$ is a gauge-type function, from Lemma \ref{theoremleft} it follows that,
  for every  $(t_0,\gamma^0_{t_0},\eta^0_{t_0})\in [\tilde{t},T]\times \Lambda^{\tilde{t}}\times \Lambda^{\tilde{t}}$ satisfy
$$
\Psi(\gamma^0_{t_0},\eta^0_{t_0})\geq \sup_{(s,\gamma_s,\eta_s)\in [\tilde{t},T]\times \Lambda^{\tilde{t}}\times \Lambda^{\tilde{t}}}\Psi(\gamma_s,\eta_s)-\frac{1}{\beta},\
\    \mbox{and} \ \ \Psi(\gamma^0_{t_0},\eta^0_{t_0})\geq \Psi(\tilde{\gamma}_{\tilde{t}},\tilde{\gamma}_{\tilde{t}}) >\frac{\tilde{m}}{2},
 $$
  there exist $(\hat{t},\hat{\gamma}_{\hat{t}},\hat{\eta}_{\hat{t}})\in [\tilde{t},T]\times \Lambda^{\tilde{t}}\times \Lambda^{\tilde{t}}$ and a sequence $\{(t_i,\gamma^i_{t_i},\eta^i_{t_i})\}_{i\geq1}\subset
  [\tilde{t},T]\times \Lambda^{\tilde{t}}\times \Lambda^{\tilde{t}}$ such that
  \begin{description}
        \item{(i)} $\Upsilon^2(\gamma^0_{t_0},\hat{\gamma}_{\hat{t}})+\Upsilon^2(\eta^0_{t_0},\hat{\eta}_{\hat{t}})+|\hat{t}-t_0|^2\leq \frac{1}{\beta}$,
         $\Upsilon^2(\gamma^i_{t_i},\hat{\gamma}_{\hat{t}})+\Upsilon^2(\eta^i_{t_i},\hat{\eta}_{\hat{t}})+|\hat{t}-t_i|^2
         \leq \frac{1}{\beta2^i}$ and $t_i\uparrow \hat{t}$ as $i\rightarrow\infty$,
        \item{(ii)}  $\Psi(\hat{\gamma}_{\hat{t}},\hat{\eta}_{\hat{t}})%-\beta[\Upsilon^2(\tilde{\gamma}_{\tilde{t}},\hat{\gamma}_{\hat{t}})+\Upsilon^2(\tilde{\gamma}_{\tilde{t}},\hat{\eta}_{\hat{t}})]
            -\sum_{i=0}^{\infty}\frac{1}{2^i}[\Upsilon^2(\gamma^i_{t_i},\hat{\gamma}_{\hat{t}})
        +\Upsilon^2(\eta^i_{t_i},\hat{\eta}_{\hat{t}})+|\hat{t}-t_i|^2]\geq \Psi(\gamma^0_{t_0},\eta^0_{t_0})$, and
        \item{(iii)}    for all $(s,\gamma_s,\eta_s)\in [\hat{t},T]\times \Lambda^{\hat{t}}\times \Lambda^{\hat{t}}\setminus \{(\hat{t},\hat{\gamma}_{\hat{t}},\hat{\eta}_{\hat{t}})\}$,
        \begin{eqnarray}\label{iii4}
        \Psi^1(\gamma_s,\eta_s)%-\beta[\Upsilon^2(\tilde{\gamma}_{\tilde{t}},\gamma_s)+\Upsilon^2(\tilde{\gamma}_{\tilde{t}},\eta_s)]
       <\Psi^1(\hat{\gamma}_{\hat{t}},\hat{\eta}_{\hat{t}}).%-\beta[\Upsilon^2(\tilde{\gamma}_{\tilde{t}},\hat{\gamma}_{\hat{t}})
            %+\Upsilon^2(\tilde{\gamma}_{\tilde{t}},\hat{\eta}_{\hat{t}})]
            %-\sum_{i=0}^{\infty}\frac{1}{2^i}[\Upsilon^2(\gamma^i_{t_i},\hat{\gamma}_{\hat{t}})
%            +\Upsilon^2(\eta^i_{t_i},\hat{\eta}_{\hat{t}})+|\hat{t}-t_i|^2].
        \end{eqnarray}
      %  \begin{eqnarray}\label{iii4}
%        &&\Psi(\gamma_s,\eta_s)%-\beta[\Upsilon^2(\tilde{\gamma}_{\tilde{t}},\gamma_s)+\Upsilon^2(\tilde{\gamma}_{\tilde{t}},\eta_s)]
%        -\sum_{i=0}^{\infty}
%        \frac{1}{2^i}[\Upsilon^2(\gamma^i_{t_i},\gamma_s)+\Upsilon^2(\eta^i_{t_i},\eta_s)+|{s}-t_i|^2]\nonumber\\
%           &<&\Psi(\hat{\gamma}_{\hat{t}},\hat{\eta}_{\hat{t}})%-\beta[\Upsilon^2(\tilde{\gamma}_{\tilde{t}},\hat{\gamma}_{\hat{t}})
%            %+\Upsilon^2(\tilde{\gamma}_{\tilde{t}},\hat{\eta}_{\hat{t}})]
%            -\sum_{i=0}^{\infty}\frac{1}{2^i}[\Upsilon^2(\gamma^i_{t_i},\hat{\gamma}_{\hat{t}})
%            +\Upsilon^2(\eta^i_{t_i},\hat{\eta}_{\hat{t}})+|\hat{t}-t_i|^2].
%        \end{eqnarray}
        where we define
$$
     \Psi^1(\gamma_t,\eta_t):=  \Psi(\gamma_t,\eta_t)
        -\sum_{i=0}^{\infty}
        \frac{1}{2^i}[{\Upsilon}^2(\gamma^i_{t_i},\gamma_t)+{\Upsilon}^2(\eta^i_{t_i},\eta_t)+|{t}-t_i|^2], \  (t,\gamma_t,\eta_t)\in [\tilde{t},T]\times  \Lambda^{\tilde{t}}\times \Lambda^{\tilde{t}}.
$$
        \end{description}
             We should note that the point
             $({\hat{t}},\hat{{\gamma}}_{{\hat{t}}},\hat{{\eta}}_{{\hat{t}}})$ depends on $\beta$ and
              $\varepsilon$.
\par
$Step\ 2.$
There exists ${{M}_0}>0$
    such that
                   \begin{eqnarray}\label{5.10jiajiaaaa}||\hat{\gamma}_{\hat{t}}||_0\vee||\hat{\eta}_{\hat{t}}||_0<M_0,
                   \end{eqnarray} % the following result  holds true:
 \begin{eqnarray}\label{5.10}
                      \beta ||\hat{{\gamma}}_{{\hat{t}}}-\hat{{\eta}}_{{\hat{t}}}||_{0}^2
                         %+\beta|\hat{{\gamma}}_{{\hat{t}}}(\hat{t})-\hat{{\eta}}_{{\hat{t}}}(\hat{t})|_m^4
                         %||\hat{{\gamma}}_{{\hat{t}}}-\hat{{\eta}}_{{\hat{t}}}||_{0}^2
                         \rightarrow0 \ \mbox{as} \ \beta\rightarrow\infty.
 \end{eqnarray}
  Let us show the above. First,   noting $\nu$ is independent of  $\beta$, by the definition of  ${\Psi}$,
 there exists an ${M}_0>0$  that is sufficiently  large   that
           $
           \Psi(\gamma_t, \eta_t)<0
           $ for all $t\in [T-\bar{a},T]$ and $||\gamma_t||_0\vee||\eta_t||_0\geq M_0$. Thus, we have $||\hat{\gamma}_{\hat{t}}||_0\vee||\hat{\eta}_{\hat{t}}||_0\vee
           ||{\gamma}^{0}_{t_{0}}||_0\vee||{\eta}^{0}_{t_{0}}||_0<M_0$.
           \par
   Second, by (\ref{iii4}), we have
    %$\Psi(\gamma_s,\eta_s)-\sum_{i=0}^{\infty}[\Upsilon^2(\gamma^i_{t_i},\gamma_s)+\Upsilon^2(\eta^i_{t_i},\eta_s)]
%            <\Psi(\hat{\gamma}_{\hat{t}},\hat{\eta}_{\hat{t}})-\sum_{i=0}^{\infty}[\Upsilon^2(\gamma^i_{t_i},\hat{\gamma}_{\hat{t}})
%            +\Upsilon^2(\eta^i_{t_i},\hat{\eta}_{\hat{t}})]$
   % the definition of $(\hat{t},\hat{\gamma}_{\hat{t}},\hat{\eta}_{\hat{t}})$, we have
 \begin{eqnarray}\label{5.56789}
                        &&2\Psi(\hat{\gamma}_{\hat{t}},\hat{\eta}_{\hat{t}})%-2\beta[\Upsilon^2(\tilde{\gamma}_{\tilde{t}},\hat{\gamma}_{\hat{t}})
            %+\Upsilon^2(\tilde{\gamma}_{\tilde{t}},\hat{\eta}_{\hat{t}})]
            -2\sum_{i=0}^{\infty}\frac{1}{2^i}[\Upsilon^2(\gamma^i_{t_i},\hat{\gamma}_{\hat{t}})
            +\Upsilon^2(\eta^i_{t_i},\hat{\eta}_{\hat{t}})+|\hat{t}-t_i|^2]\nonumber\\
            &\geq&  \Psi(\hat{\gamma}_{\hat{t}},\hat{\gamma}_{\hat{t}})%-\beta[\Upsilon^2(\tilde{\gamma}_{\tilde{t}},\hat{\gamma}_{\hat{t}})
            %+\Upsilon^2(\tilde{\gamma}_{\tilde{t}},\hat{\gamma}_{\hat{t}})]
            -\sum_{i=0}^{\infty}\frac{1}{2^i}[\Upsilon^2(\gamma^i_{t_i},\hat{\gamma}_{\hat{t}})
            +\Upsilon^2(\eta^i_{t_i},\hat{\gamma}_{\hat{t}})+|\hat{t}-t_i|^2]\nonumber\\
            &&
            +\Psi(\hat{\eta}_{\hat{t}},\hat{\eta}_{\hat{t}})%-\beta[\Upsilon^2(\tilde{\eta}_{\tilde{t}},\hat{\gamma}_{\hat{t}})
            %+\Upsilon^2(\tilde{\gamma}_{\tilde{t}},\hat{\eta}_{\hat{t}})]
            -\sum_{i=0}^{\infty}\frac{1}{2^i}[\Upsilon^2(\gamma^i_{t_i},\hat{\eta}_{\hat{t}})
            +\Upsilon^2(\eta^i_{t_i},\hat{\eta}_{\hat{t}})+|\hat{t}-t_i|^2].
 \end{eqnarray}
 This implies that
 \begin{eqnarray}\label{5.6}
                         2{\beta}\Upsilon^2(\hat{\gamma}_{\hat{t}},\hat{\eta}_{\hat{t}})
                         \leq|W_1(\hat{\gamma}_{\hat{t}})-W_1(\hat{\eta}_{\hat{t}})|
                                   +|W_2(\hat{\gamma}_{\hat{t}})-W_2(\hat{\eta}_{\hat{t}})|+
                                   \sum_{i=0}^{\infty}\frac{1}{2^i}[\Upsilon^2(\eta^i_{t_i},\hat{\gamma}_{\hat{t}})+\Upsilon^2(\gamma^i_{t_i},\hat{\eta}_{\hat{t}})].
                                  % &\leq& 2L(2+||\hat{{\gamma}}^{1}_{{\hat{t}}}||_0+||\hat{{\gamma}}_{{\hat{t}}}^{2}||_0)
%                                   \leq 4L(1+M^m).
 \end{eqnarray}
 On the other hand, by Lemma \ref{theoremS00044},
 \begin{eqnarray}\label{4.7jiajia130}
 &&\sum_{i=0}^{\infty}\frac{1}{2^i}[\Upsilon^2(\eta^i_{t_i},\hat{\gamma}_{\hat{t}})+\Upsilon^2(\gamma^i_{t_i},\hat{\eta}_{\hat{t}})]\nonumber\\
 %\leq
% 2\sum_{i=0}^{\infty}\frac{1}{2^i}[||\eta^i_{t_i}-\hat{\gamma}_{\hat{t}}||_0^4+||\gamma^i_{t_i}-\hat{\eta}_{\hat{t}}||_0^4]\\
 %&\leq&2^4\sum_{i=0}^{\infty}\frac{1}{2^i}[||\eta^i_{t_i}-\hat{\eta}_{\hat{t}}||_0^4%+||\hat{\eta}_{\hat{t}}-\hat{\gamma}_{\hat{t}}||_0^4
% +||\gamma^i_{t_i}-\hat{\gamma}_{\hat{t}}||_0^4+2||\hat{\gamma}_{\hat{t}}-\hat{\eta}_{\hat{t}}||_0^4]\\
 %&\leq&2^4\sum_{i=0}^{\infty}\frac{1}{2^i}[||\eta^i_{t_i}-\hat{\eta}_{\hat{t}}||_0^4%+||\hat{\eta}_{\hat{t}}-\hat{\gamma}_{\hat{t}}||_0^4
% +||\gamma^i_{t_i}-\hat{\gamma}_{\hat{t}}||_0^4+2^4||{\gamma}^{0}_{t_{0}}-\tilde{\gamma}_{\tilde{t}}||_0^4+2^7||{\eta}^{0}_{t_{0}}-\hat{\eta}_{\hat{t}}||_0^4
% +2^7||{\gamma}^{0}_{t_{0}}-{\eta}^{0}_{t_{0}}||_0^4]\\
 &\leq&2\sum_{i=0}^{\infty}\frac{1}{2^i}[\Upsilon^2(\eta^i_{t_i},\hat{\eta}_{\hat{t}})%+||\hat{\eta}_{\hat{t}}-\hat{\gamma}_{\hat{t}}||_0^4
 +\Upsilon^2(\gamma^i_{t_i},\hat{\gamma}_{\hat{t}})+2\Upsilon^2(\hat{\gamma}_{\hat{t}},\hat{\eta}_{\hat{t}})]
 \leq\frac{4}{\beta}+{8}\Upsilon^2(\hat{\gamma}_{\hat{t}},\hat{\eta}_{\hat{t}}).
 %\frac{54M_1}{\beta}(2^4+\frac{1}{3})=\frac{882M_1}{\beta}.
 \end{eqnarray}
Thus we have
 \begin{eqnarray}\label{5.jia6}
                         &&(2{\beta}-8)\Upsilon^2(\hat{\gamma}_{\hat{t}},\hat{\eta}_{\hat{t}})
                         \leq |W_1(\hat{\gamma}_{\hat{t}})-W_1(\hat{\eta}_{\hat{t}})|
                                   +|W_2(\hat{\gamma}_{\hat{t}})-W_2(\hat{\eta}_{\hat{t}})|+\frac{4}{\beta}\nonumber\\
                                   %&\leq&|W_1(\hat{\gamma}_{\hat{t}})-W_1(\hat{\eta}_{\hat{t}})|
%                                   +|W_2(\hat{\gamma}_{\hat{t}})-W_2(\hat{\eta}_{\hat{t}})|
%                                   + 54\times 2^7\sum_{i=0}^{\infty}\frac{1}{2^i}[\frac{1}{2^i\beta}+\frac{1}{\beta}+16M_0^4]\nonumber\\
                                   &\leq& 2L(2+||\hat{\gamma}_{\hat{t}}||_0+||\hat{\eta}_{\hat{t}}||_0)+\frac{4}{\beta}
                                   \leq 4L(1+M_0)+\frac{4}{\beta}.
 \end{eqnarray}
   Letting $\beta\rightarrow\infty$, we get
                $$\Upsilon^2(\hat{\gamma}_{\hat{t}},\hat{\eta}_{\hat{t}})
                \rightarrow0\
                          \mbox{as} \ \beta\rightarrow+\infty.
                          $$
 Then from (\ref{s0}) and (\ref{5.jia6}) it follows that (\ref{5.10}) holds.
 \par
   $Step\ 3.$ There exists $N_0>0$ such that
                   $\hat{t}\in [T-\bar{a},T)$ for all $\beta\geq N_0$.
 \par
By (\ref{5.10}), we can let $N_0>0$ be a large number such that
$$
                         L||\hat{\gamma}_{\hat{t}}-\hat{\eta}_{\hat{t}}||_0
                         \leq
                         \frac{\tilde{m}}{4},
$$
               for all $\beta\geq N_0$.
            Then we have $\hat{t}\in [T-\bar{a},T)$ for all $\beta\geq N_0$. Indeed, if say $\hat{t}=T$,  we will deduce the following contradiction:
 \begin{eqnarray*}
                         \frac{\tilde{m}}{2}\leq\Psi(\hat{\gamma}_{\hat{t}},\hat{\eta}_{\hat{t}})\leq \phi(\hat{\gamma}_{\hat{t}})-\phi(\hat{\eta}_{\hat{t}})\leq
                        L||\hat{\gamma}_{\hat{t}}-\hat{\eta}_{\hat{t}}||_0
                         \leq
                         \frac{\tilde{m}}{4}.
 \end{eqnarray*}
%
%
%
%
%\par
%Next, noting $\Psi(\gamma^1_T,\gamma^2_T)\leq0<\frac{\tilde{m}}{2}\leq \Psi(\tilde{\gamma}_{\tilde{t}},\tilde{\gamma}_{\tilde{t}})$ for every $\gamma^1_T,\gamma^2_T\in \Lambda_T$,
%   we have $\hat{t}\in (T-\bar{a},T)$.
%\par
%Next, from (\ref{psi0}), it follows that
% \begin{eqnarray*}
% \frac{\tilde{m}}{2}&\leq& W(\hat{{\gamma}}^{1}_{{\hat{t}}})- V(\hat{{\gamma}}^{2}_{{\hat{t}}})\\
% &=&W(\hat{{\gamma}}^{1}_{{\hat{t}}})-W(\hat{{\gamma}}^{1}_{{\hat{t}},T})+W(\hat{{\gamma}}^{1}_{{\hat{t}},T})
%    -V(\hat{{\gamma}}^{2}_{{\hat{t}},T})+V(\hat{{\gamma}}^{2}_{{\hat{t}},T})- V(\hat{{\gamma}}^{2}_{{\hat{t}}})\\
%    &\leq&(C+C_1)(1+2M_0)(T-\hat{t})^{\frac{1}{2}}+L\mu(|\hat{{\gamma}}^{1}_{{\hat{t}}}-\hat{{\gamma}}^{2}_{{\hat{t}}}|).
% \end{eqnarray*}
%On the other hand, from (\ref{5.10}) it follows that there exists a constant $N>0$ such that
% $$
%  ||\hat{{\gamma}}^{1}_{{\hat{t}}}-\hat{{\gamma}}^{2}_{{\hat{t}}}||_0<\frac{\tilde{m}}{4L(1+\mu([0,T]))}, \ \ \mbox{for  all} \ \ \beta>N.
% $$
% Then
% $$
% (C+C_1)(1+2M_0)(T-\hat{t})^{\frac{1}{2}}>\frac{\tilde{m}}{4}.
% $$
% %Letting $\Delta=\frac{\tilde{m}^2}{16(C+C_1)^2(1+2M_0)^2}$,
%  Thus, we have
%  $\hat{t} \in [T-\bar{a},T)$  for all $\beta>N$.
\par
 $Step\ 4.$    Completion of the proof.
\par
          From above all,  for the fixed   $N_0>0$ in step 3, %and $0<\Delta<T$ in step 3,
          we  find
$(\hat{t},\hat{\gamma}_{\hat{t}}), (\hat{t},\hat{\eta}_{\hat{t}})\in [\tilde{t}, T]\times
                 \Lambda^{\tilde{t}}$   satisfying $\tilde{t}\in [T-\bar{a},T)$  for all $\beta\geq N_0$
           such that
\begin{eqnarray}\label{psi4}
            \Psi^1(\hat{\gamma}_{\hat{t}},\hat{\eta}_{\hat{t}})\geq \Psi(\tilde{\gamma}_{\tilde{t}},\tilde{\gamma}_{\tilde{t}}) \ \ \mbox{and} \ \      \Psi^1(\hat{\gamma}_{\hat{t}},\hat{\eta}_{\hat{t}})\geq
                   {\Psi}^1(\gamma_t,\eta_t),
                   \  (t,\gamma_t,\eta_t)\in [\hat{t},T]\times  \Lambda^{\hat{t}}\times \Lambda^{\hat{t}}.
\end{eqnarray}
%where we define
%$$
%     \Psi^1(\gamma_t,\eta_t):=  \Psi(\gamma_t,\eta_t)
%        -\sum_{i=0}^{\infty}
%        \frac{1}{2^i}[{\Upsilon}^2(\gamma^i_{t_i},\gamma_t)+{\Upsilon}^2(\eta^i_{t_i},\eta_t)+|{t}-t_i|^2], \ \  \  (t,\gamma_t,\eta_t)\in [\tilde{t},T]\times  \Lambda^{\tilde{t}}\times \Lambda^{\tilde{t}}.
%$$
%We put, for $(t,\gamma_t,\eta_t)\in [T-\bar{a},T]\times {\Lambda}\times {\Lambda}$,
 Now we consider the function,
              for $(t,\gamma_t), (s,\eta_s)\in [\hat{t},T]\times{{\Lambda}}$,
\begin{eqnarray}\label{4.1111}
                 \Psi_{\delta}(\gamma_t,\eta_s)=W'_{1}(\gamma_t)-W'_{2}(\eta_s)-2\beta(\Upsilon^2(\gamma_{t},\hat{{\xi}}_{{\hat{t}}})+\Upsilon^2(\eta_{s},\hat{{\xi}}_{{\hat{t}}}))-\frac{1}{\delta}|s-t|^2,
\end{eqnarray}
where
\begin{eqnarray*}
                             {W}'_{1}(\gamma_t)=W_1(\gamma_t)-\varepsilon\frac{\nu T-t}{\nu
                 T}\Upsilon^2(\gamma_t)
                % -\varepsilon |t-{\hat{t}}|^2
                 -\varepsilon \overline{\Upsilon}^2(\gamma_t,\hat{\gamma}_{\hat{t}})-\sum_{i=0}^{\infty}
        \frac{1}{2^i}\overline{\Upsilon}^2(\gamma^i_{t_i},\gamma_t),
        \end{eqnarray*}
        \begin{eqnarray*}
                             {W}'_{2}(\eta_s)&=&W_2(\eta_s)+\varepsilon\frac{\nu T-s}{\nu
                 T}\Upsilon^2(\eta_s)
                 +\varepsilon \overline{\Upsilon}^2(\eta_s,\hat{\eta}_{\hat{t}})+\sum_{i=0}^{\infty}
        \frac{1}{2^i}{\Upsilon}^2(\eta^i_{t_i},\eta_s),
\end{eqnarray*}
and
$$
\hat{\xi}_{\hat{t}}=\frac{\hat{\gamma}_{\hat{t}}+\hat{\eta}_{\hat{t}}}{2}.
$$
 Define a sequence of positive numbers $\{\delta_i\}_{i\geq0}$  by
        $\delta_i=\frac{1}{2^i}$ for all $i\geq0$.  For every $\delta>0$,
           from Lemma \ref{theoremleft} it follows that,
  for every  $(\check{t}_{0},\check{\gamma}^{0}_{\check{t}_{0}}), (\check{s}_{0},\check{\eta}^{0}_{\check{s}_{0}})\in [\hat{t},T]\times  \Lambda^{\hat{t}}$ satisfy
$$
\Psi_{\delta}(\check{\gamma}^{0}_{\check{t}_{0}},\check{\eta}^{0}_{\check{s}_{0}})\geq \sup_{(t,\gamma_t),(s,\eta_s)\in [\hat{t},T]\times \Lambda^{\hat{t}}}\Psi_{\delta}(\gamma_t,\eta_s)-{\delta},\
%\    \mbox{and} \ \ \Psi_{\delta}(\gamma^0_{t_0},\eta^0_{t_0})\geq \Psi(\tilde{\gamma}_{\tilde{t}},\tilde{\gamma}_{\tilde{t}}) >\frac{\tilde{m}}{2},
 $$
  there exist $(\check{t},\check{\gamma}_{\check{t}}), (\check{s},\check{\eta}_{\check{s}})\in [\hat{t},T]\times \Lambda^{\hat{t}}$ and two sequences $\{(\check{t}_{i},\check{\gamma}^{i}_{\check{t}_{i}})\}_{i\geq1}, \{(\check{s}_{i},\check{\eta}^{i}_{\check{s}_{i}})\}_{i\geq1}\subset
  [\hat{t},T]\times \Lambda^{\hat{t}}$ such that
  \begin{description}
        \item{(i)} $\overline{\Upsilon}^2(\check{\gamma}^{0}_{\check{t}_{0}},\check{\gamma}_{\check{t}})+\overline{\Upsilon}^2(\check{\eta}^0_{\check{s}_0},\check{\eta}_{\check{s}})\leq {\delta}$,
         $\overline{\Upsilon}^2(\check{\gamma}^{i}_{\check{t}_{i}},\check{\gamma}_{\check{t}})
         +\overline{\Upsilon}^2(\check{\eta}^{i}_{\check{s}_{i}},\check{\eta}_{\check{s}})\leq \frac{\delta}{2^i}$
          and $\check{t}_{i}\uparrow \check{t}$, $\check{s}_{i}\uparrow \check{s}$ as $i\rightarrow\infty$,
        \item{(ii)}  $\Psi_{\delta}(\check{\gamma}_{\check{t}},\check{\eta}_{\check{s}})%-\beta[\overline{\Upsilon}^2(\tilde{\gamma}_{\tilde{t}},\hat{\gamma}^k_{\hat{t}_k})+\overline{\Upsilon}^2(\tilde{\gamma}_{\tilde{t}},\hat{\eta}^k_{\hat{t}_k})]
            -\sum_{i=0}^{\infty}\frac{1}{2^i}[\overline{\Upsilon}^2(\check{\gamma}^{i}_{\check{t}_{i}},\check{\gamma}_{\check{t}})
        +\overline{\Upsilon}^2(\check{\eta}^{i}_{\check{s}_{i}},\check{\eta}_{\check{s}})]\geq \Psi_{\delta}(\check{\gamma}^{0}_{\check{t}_{0}},\check{\eta}^{0}_{\check{s}_{0}})$, and
        \item{(iii)}    for all $(t,\gamma_t,s,\eta_s)\in [\check{t},T]\times \Lambda^{\check{t}}\times[\check{s},T]\times \Lambda^{\check{s}}\setminus \{(\check{t},\check{\gamma}_{\check{t}},\check{s},\check{\eta}_{\check{s}})\}$,
        \begin{eqnarray}\label{519jia}
        &&\Psi_{\delta}(\gamma_t,\eta_s)%-\beta[\overline{\Upsilon}^2(\tilde{\gamma}_{\tilde{t}},\gamma_s)+\overline{\Upsilon}^2(\tilde{\gamma}_{\tilde{t}},\eta_s)]
        -\sum_{i=0}^{\infty}
        \frac{1}{2^i}[\overline{\Upsilon}^2(\check{\gamma}^{i}_{\check{t}_{i}},\gamma_t)+\overline{\Upsilon}^2(\check{\eta}^{i}_{\check{s}_{i}},\eta_s)]\nonumber\\
            &<&\Psi_{\delta}(\check{\gamma}_{\check{t}},\check{\eta}_{\check{s}})%-\beta[\overline{\Upsilon}^2(\tilde{\gamma}_{\tilde{t}},\hat{\gamma}^k_{\hat{t}_k})
            %+\overline{\Upsilon}^2(\tilde{\gamma}_{\tilde{t}},\hat{\eta}^k_{\hat{t}_k})]
            -\sum_{i=0}^{\infty}\frac{1}{2^i}[\overline{\Upsilon}^2(\check{\gamma}^{i}_{\check{t}_{i}},\check{\gamma}_{\check{t}})
            +\overline{\Upsilon}^2(\check{\eta}^{i}_{\check{s}_{i}},\check{\eta}_{\check{s}})].
        \end{eqnarray}
        \end{description}
  By the following Lemma \ref{lemma4.4}, we have
\begin{eqnarray}\label{4.23}
\lim_{\delta\rightarrow0}[\overline{\Upsilon}^2(\check{\gamma}_{\check{t}},\hat{\gamma}_{\hat{t}})
+\overline{\Upsilon}^2(\check{\eta}_{\check{s}},\hat{\eta}_{\hat{t}})]=0.
\end{eqnarray}
From (\ref{4.23}) and ${\hat{t}}<T$
 for $\beta>N_0$, it follows that, for every fixed $\beta>N_0$,   constant $ K_\beta>0$ exists such that
$$
            |\check{t}|\vee|\check{s}|<T, %\ \    |x^{k}_0|\vee|y^{k}_0|<\frac{M_0}{2},
            \ \ \mbox{for all}    \ \ 0<\delta< K_\beta.
$$
Now, for every $\beta>N_0$ and $0<\delta< K_\beta$, from the definition of viscosity solutions it follows that
\begin{eqnarray}\label{vis1}
                      &&
                      -\frac{\varepsilon}{\nu T}\Upsilon^2(\check{\gamma}_{{\check{t}}})
                     +2\varepsilon({\check{t}}-{\hat{t}})
                     +2\sum_{i=0}^{\infty}\frac{1}{2^i}[(\check{t}-\check{t}_{i})+(\check{t}-t_i)]+\frac{2}{\delta}(\check{t}-\check{s})\nonumber\\
                     &&+{\mathbf{H}}(\check{\gamma}_{{\check{t}}}, 2\beta\partial_x\Upsilon^2(\check{\gamma}_{{\check{t}}}-\hat{\xi}_{{\hat{t}},\check{t},A})
                      +\varepsilon\partial_x\Upsilon^2(\check{\gamma}_{{\check{t}}}-\hat{\gamma}_{{\hat{t}},\check{t},A})
                      +\varepsilon\frac{\nu T-{\check{t}}}{\nu T}\partial_x\Upsilon^2(\check{\gamma}_{{\check{t}}})\nonumber\\
                     && +\sum_{i=0}^{\infty}\frac{1}{2^i}
                      \partial_x\Upsilon^2(\check{\gamma}_{\check{t}}-\check{\gamma}^{i}_{\check{t}_{i},\check{t},A})
                      +\sum_{i=0}^{\infty}\frac{1}{2^i}\partial_x\Upsilon^2(\check{\gamma}_{\check{t}}-\gamma^i_{t_i,\check{t},A})
                                    )\geq c;
 \end{eqnarray}
 and
 \begin{eqnarray}\label{vis2}
                     &&
                     \frac{\varepsilon}{\nu T}\Upsilon^2(\check{\eta}_{{\check{s}}})
                     -2\varepsilon({\check{s}}-{\hat{t}}) -2\sum_{i=0}^{\infty}\frac{1}{2^i}(\check{s}-\check{s}_{i})+\frac{2}{\delta}(\check{t}-\check{s})\nonumber\\
                      &&
                     +{\mathbf{H}}(\check{\eta}_{{\check{s}}},
                      -2\beta\partial_x\Upsilon^2(\check{\eta}_{{\check{s}}}-\hat{\xi}_{{\hat{t}},\check{s},A})
                      -\varepsilon\partial_x\Upsilon^2(\check{\eta}_{{\check{s}}}-\hat{\eta}_{{\hat{t}},\check{s},A})
                       -\varepsilon\frac{\nu T-{\check{s}}}{\nu T}
                     \partial_x\Upsilon^2(\check{\eta}_{{\check{s}}})\nonumber\\
                      &&-\sum_{i=0}^{\infty}\frac{1}{2^i}
                      \partial_x\Upsilon^2(\check{\eta}_{\check{s}}-\check{\eta}^{i}_{\check{s}_{i},\check{s},A})
                      -\sum_{i=0}^{\infty}\frac{1}{2^i}\partial_x\Upsilon^2(\check{\eta}_{\check{s}}-\eta^i_{t_i,\check{s},A})
                                     )\leq0.
  \end{eqnarray}
 % Here and in the sequel,
% for notational simplicity,
%we use $\partial_x\Upsilon^2(\cdot,\cdot)$  to denote the first spatial derivative with respect to the first variable.
  We notice that, by  the property (i) of $(\check{t},\check{\gamma}_{\check{t}},\check{s},\check{\eta}_{\check{s}})$,
  \begin{eqnarray*}
  2\sum_{i=0}^{\infty}\frac{1}{2^i}[(\check{s}-\check{s}_{i})+(\check{t}-\check{t}_{i})]
  \leq4\sum_{i=0}^{\infty}\frac{1}{2^i}\bigg{(}\frac{\delta}{2^i}\bigg{)}^{\frac{1}{2}}\leq 8\delta^{\frac{1}{2}};
    \end{eqnarray*}
    \begin{eqnarray*}
    |\partial_x\Upsilon^2(\check{\gamma}_{{\check{t}}}-\hat{\gamma}_{{\hat{t}},\check{t},A})|
    +|\partial_x\Upsilon^2(\check{\eta}_{{\check{s}}}-\hat{\gamma}_{{\hat{t}},\check{s},A})|\leq 4|e^{(\check{t}-\hat{t})A}\hat{\gamma}_{{\hat{t}}}(\hat{t})-\check{\gamma}_{\check{t}}(\check{t})|
                      +4|e^{(\check{s}-\hat{t})A}\hat{\gamma}_{{\hat{t}}}(\hat{t})-\check{\eta}_{\check{s}}(\check{s})|;
    \end{eqnarray*}
    and
     \begin{eqnarray*}
  &&\bigg{|}\sum_{i=0}^{\infty}\frac{1}{2^i}
                      \partial_x\Upsilon^2(\check{\gamma}_{\check{t}}-\check{\gamma}^{i}_{\check{t}_{i},\check{t},A})\bigg{|}
                      +\bigg{|}\sum_{i=0}^{\infty}\frac{1}{2^i}
                      \partial_x\Upsilon^2(\check{\eta}_{\check{s}}-\check{\eta}^{i}_{\check{s}_{i},\check{s},A})\bigg{|}\\
                      &\leq&4\sum_{i=0}^{\infty}\frac{1}{2^i}[|e^{(\check{t}-\check{t}_{i})A}\check{\gamma}^{i}_{\check{t}_{i}}(\check{t}_{i})-\check{\gamma}_{\check{t}}(\check{t})|
                      +|e^{(\check{s}-\check{s}_{i})A}\check{\eta}^{i}_{\check{s}_{i}}(\check{s}_{i})-\check{\eta}_{\check{s}}(\check{s})|]
                      \leq
                      8\sum_{i=0}^{\infty}\frac{1}{2^i}\bigg{(}\frac{\delta}{2^i}\bigg{)}^{\frac{1}{2}}\leq 16{\delta}^{\frac{1}{2}}.
                        \end{eqnarray*}
  Combining(\ref{vis1}) and (\ref{vis2}), and letting $\delta\rightarrow0$,   we obtain
  \begin{eqnarray}\label{vis112}
                     c+ \frac{\varepsilon}{\nu T}(\Upsilon^2(\hat{{\gamma}}_{{\hat{t}}})+\Upsilon^2(\hat{{\eta}}_{{\hat{t}}}))
                     \leq{\mathbf{H}}_1-{\mathbf{H}}_2
                     +2\sum_{i=0}^{\infty}\frac{1}{2^i}(\hat{t}-t_i),
\end{eqnarray}
  where
\begin{eqnarray*}
                       {\mathbf{H}}_1={\mathbf{H}}(\hat{{\gamma}}_{{\hat{t}}},
                                      2\beta\partial_x\Upsilon^2(\hat{\gamma}_{{\hat{t}}}-\hat{\xi}_{{\hat{t}}})
                                      +\varepsilon\frac{\nu T-{\hat{t}}}{\nu T}\partial_x\Upsilon^2(\hat{{\gamma}}_{{\hat{t}}})
                                     +\sum_{i=0}^{\infty}\frac{1}{2^i}\partial_x\Upsilon^2(\hat{\gamma}_{\hat{t}}-\gamma^i_{t_i,\hat{t},A}))
                                     ;\nonumber\\
                      {\mathbf{H}}_2={\mathbf{H}}(\hat{{\eta}}_{{\hat{t}}},
                      -2\beta\partial_x\Upsilon^2(\hat{\eta}_{{\hat{t}}}-\hat{\xi}_{{\hat{t}}})-\varepsilon\frac{\nu T-{\hat{t}}}{\nu T}\partial_x\Upsilon^2(\hat{{\eta}}_{{\hat{t}}})
                                     -\sum_{i=0}^{\infty}\frac{1}{2^i}
                                     \partial_x\Upsilon^2(\hat{\eta}_{\hat{t}}-\eta^i_{t_i,\hat{t},A})).
\end{eqnarray*}
 On the other hand, by  a simple calculation we obtain
\begin{eqnarray}\label{v4}
                {\mathbf{H}}_1-{\mathbf{H}}_2
                \leq\sup_{u\in U}(J_{1}+J_{2}),
\end{eqnarray}
 where
\begin{eqnarray}\label{j1}
                               J_{1}&=&( {F}(\hat{{\gamma}}_{{\hat{t}}},u),2\beta\partial_x\Upsilon^2(\hat{\gamma}_{{\hat{t}}}-\hat{\xi}_{{\hat{t}}})
                               +\varepsilon\frac{\nu T-{\hat{t}}}{\nu T}\partial_x\Upsilon^2(\hat{{\gamma}}_{{\hat{t}}})
                                     +\sum_{i=0}^{\infty}\frac{1}{2^i}\partial_x\Upsilon^2(\hat{\gamma}_{\hat{t}}-\gamma^i_{t_i,\hat{t},A}))_{R^d}\nonumber\\
                                     &&  -( {F}(\hat{{\eta}}_{{\hat{t}}},u),
                                            -2\beta\partial_x\Upsilon^2(\hat{\eta}_{{\hat{t}}}-\hat{\xi}_{{\hat{t}}})-\varepsilon\frac{\nu T-{\hat{t}}}{\nu T}\partial_x\Upsilon^2(\hat{{\eta}}_{{\hat{t}}})
                                     -\sum_{i=0}^{\infty}\frac{1}{2^i}\partial_x\Upsilon^2(\hat{\eta}_{\hat{t}}-\eta^i_{t_i,\hat{t},A}))_{R^d}\nonumber\\
                                  &\leq&4\beta{L}||\hat{{\gamma}}_{{\hat{t}}}-\hat{{\eta}}_{{\hat{t}}}||_0^2
                                  +8\varepsilon \frac{\nu T-{\hat{t}}}{\nu T} L(1+||\hat{{\gamma}}_{{\hat{t}}}||^2_0+||\hat{{\eta}}_{{\hat{t}}}||^2_0)
                                            \nonumber\\
                                           && +4L\sum_{i=0}^{\infty}\frac{1}{2^i}[|e^{(\hat{t}-t_i)A}\gamma^i_{t_i}(t_i)-\hat{\gamma}_{\hat{t}}(\hat{t})|
                                            +|e^{(\hat{t}-t_i)A}\eta^i_{t_i}(t_i)-\hat{\eta}_{\hat{t}}(\hat{t})|]
                                            (1+||\hat{{\gamma}}_{{\hat{t}}}||_0+||\hat{{\eta}}_{{\hat{t}}}||_0); \ \ \
\end{eqnarray}
\begin{eqnarray}\label{j3}
                                 J_{2}&=&q(\hat{{\gamma}}_{{\hat{t}}}, u)-
                                 q(\hat{{\eta}}_{{\hat{t}}},u)
                                     \leq
                                 L||\hat{{\gamma}}_{{\hat{t}}}-\hat{{\eta}}_{{\hat{t}}}||_0;
\end{eqnarray}
 We notice that, by  the property (i) of $(\hat{t},\hat{\gamma}_{\hat{t}},\hat{\eta}_{\hat{t}})$,
  \begin{eqnarray*}
2\sum_{i=0}^{\infty}\frac{1}{2^i}(\hat{t}-t_i)
  \leq2\sum_{i=0}^{\infty}\frac{1}{2^i}\bigg{(}\frac{1}{2^i\beta}\bigg{)}^{\frac{1}{2}}\leq 4{\bigg{(}\frac{1}{{\beta}}\bigg{)}}^{\frac{1}{2}},
    \end{eqnarray*}
    and
     \begin{eqnarray*}
  \sum_{i=0}^{\infty}\frac{1}{2^i}[|e^{(\hat{t}-t_i)A}\gamma^i_{t_i}(t_i)-\hat{\gamma}_{\hat{t}}(\hat{t})|
                                            +|e^{(\hat{t}-t_i)A}\eta^i_{t_i}(t_i)-\hat{\eta}_{\hat{t}}(\hat{t})|]
                      \leq2\sum_{i=0}^{\infty}\frac{1}{2^i}\bigg{(}\frac{1}{2^i\beta}\bigg{)}^{\frac{1}{2}}\leq 4\bigg{(}\frac{1}{\beta}\bigg{)}^{\frac{1}{2}}.
                        \end{eqnarray*}
 Combining (\ref{vis112})-(\ref{j3}), and letting $\beta\rightarrow\infty$, by (\ref{5.10jiajiaaaa}) and (\ref{5.10}) we obtain
 \begin{eqnarray}\label{vis122}
                     c
                                             &\leq&
                     \overline{ \lim_{\beta\rightarrow\infty}}\bigg{[}-\frac{\varepsilon}{\nu T}(\Upsilon^2(\hat{{\gamma}}_{{\hat{t}}})
                     +\Upsilon^2(\hat{{\eta}}_{{\hat{t}}}))+ 8\varepsilon \frac{\nu T-{\hat{t}}}{\nu T} L(1+||\hat{{\gamma}}_{{\hat{t}}}||^2_0+||\hat{{\eta}}_{{\hat{t}}}||^2_0)\bigg{]}.
\end{eqnarray}
Recalling $
         \nu=1+\frac{1}{16TL}
$ and $\bar{a}=\frac{1}{16L}\wedge{T}$, by  (\ref{5.3}), the following contradiction is induced:
\begin{eqnarray*}\label{vis122}
                     c\leq
                              \frac{\varepsilon}{\nu
                              T}\leq \frac{c}{2}.
\end{eqnarray*}

 The proof is now complete.
 \ \ $\Box$\par
 To complete the previous proof, it remains to state and prove the following lemma.
 \begin{lemma}\label{lemma4.4}\ \ The maximum points $(\check{t},\check{\gamma}_{\check{t}},\check{s},\check{\eta}_{\check{s}})$
of $\Psi_{\delta}(\gamma_t,\eta_s)
        -\sum_{i=0}^{\infty}
        \frac{1}{2^i}[\overline{\Upsilon}^2(\check{\gamma}^{i}_{\check{t}_{i}},\gamma_t)+\overline{\Upsilon}^2(\check{\eta}^{i}_{\check{t}_{i}},\eta_s)]$ defined by (\ref{4.1111}) in $[\check{t},T]\times \Lambda^{\check{t}}\times[\check{s},T]\times \Lambda^{\check{s}}$ satisfy  condition (\ref{4.23}).
\end{lemma}
\par
   {\bf  Proof  }. \ \
    Without loss of generality, we assume $\check{t}\leq \check{s}$.
   By (\ref{519jia}), we have
   %\begin{eqnarray}\label{519jia}
%        \Psi_{\delta}(\gamma_t,\eta_s)%-\beta[\overline{\Upsilon}^2(\tilde{\gamma}_{\tilde{t}},\gamma_s)+\overline{\Upsilon}^2(\tilde{\gamma}_{\tilde{t}},\eta_s)]
%        -\sum_{i=0}^{\infty}
%        \frac{1}{2^i}[\overline{\Upsilon}^2(\check{\gamma}^{i}_{\check{t}_{i}},\gamma_t)+\overline{\Upsilon}^2(\check{\eta}^{i}_{\check{s}_{i}},\eta_s)]
%            <\Psi_{\delta}(\check{\gamma}_{\check{t}},\check{\eta}_{\check{s}})%-\beta[\overline{\Upsilon}^2(\tilde{\gamma}_{\tilde{t}},\hat{\gamma}^k_{\hat{t}_k})
%            %+\overline{\Upsilon}^2(\tilde{\gamma}_{\tilde{t}},\hat{\eta}^k_{\hat{t}_k})]
%            -\sum_{i=0}^{\infty}\frac{1}{2^i}[\overline{\Upsilon}^2(\check{\gamma}^{i}_{\check{t}_{i}},\check{\gamma}_{\check{t}})
%            +\overline{\Upsilon}^2(\check{\eta}^{i}_{\check{s}_{i}},\check{\eta}_{\check{s}})].
%        \end{eqnarray}
    %$\Psi(\gamma_s,\eta_s)-\sum_{i=0}^{\infty}[\Upsilon^2(\gamma^i_{t_i},\gamma_s)+\Upsilon^2(\eta^i_{t_i},\eta_s)]
%            <\Psi(\hat{\gamma}_{\hat{t}},\hat{\eta}_{\hat{t}})-\sum_{i=0}^{\infty}[\Upsilon^2(\gamma^i_{t_i},\hat{\gamma}_{\hat{t}})
%            +\Upsilon^2(\eta^i_{t_i},\hat{\eta}_{\hat{t}})]$
   % the definition of $(\hat{t},\hat{\gamma}_{\hat{t}},\hat{\eta}_{\hat{t}})$, we have
 \begin{eqnarray}\label{5.56789}
                        &&2\Psi_{\delta}(\check{\gamma}_{\check{t}},\check{\eta}_{\check{s}})
            -2\sum_{i=0}^{\infty}\frac{1}{2^i}[\overline{\Upsilon}^2(\check{\gamma}^{i}_{\check{t}_{i}},\check{\gamma}_{\check{t}})
            +\overline{\Upsilon}^2(\check{\eta}^{i}_{\check{s}_{i}},\check{\eta}_{\check{s}})]\nonumber\\
            &\geq&  \Psi_{\delta}(\check{\gamma}_{\check{t},\check{s},A},\check{\gamma}_{\check{t},\check{s},A})+\Psi_{\delta}(\check{\eta}_{\check{s}},\check{\eta}_{\check{s}})
            -\sum_{i=0}^{\infty}\frac{1}{2^i}[\overline{\Upsilon}^2(\check{\gamma}^{i}_{\check{t}_{i}},\check{\gamma}_{\check{t},\check{s},A})
            +\overline{\Upsilon}^2(\check{\eta}^{i}_{\check{s}_{i}},\check{\gamma}_{\check{t},\check{s},A})\nonumber\\
            &&
            +\overline{\Upsilon}^2(\check{\gamma}^{i}_{\check{t}_{i}},\check{\eta}_{\check{s}})+\overline{\Upsilon}^2(\check{\eta}^{i}_{\check{s}_{i}},\check{\eta}_{\check{s}})]. \ \
 \end{eqnarray}
 This implies that
 \begin{eqnarray}\label{5.6}
                         \frac{2}{\delta}|\check{t}-\check{s}|^2
                         \leq|W'_1(\check{\gamma}_{\check{t},\check{s},A})-W'_1(\check{\eta}_{\check{s}})|
                                   +|W'_2(\check{\gamma}_{\check{t},\check{s},A})-W'_2(\check{\eta}_{\check{s}})|+
                                   \sum_{i=0}^{\infty}\frac{1}{2^i}[\overline{\Upsilon}^2(\check{\eta}^{i}_{\check{s}_{i}},\check{\gamma}_{\check{t},\check{s},A})
                                   +\overline{\Upsilon}^2(\check{\gamma}^{i}_{\check{t}_{i}},\check{\eta}_{\check{s}})
           ]. \ \ \
 \end{eqnarray}
 Letting $\delta\rightarrow0$, we have
 $$
 |\check{t}-\check{s}|\rightarrow0 \ \ \mbox{as}\  \ \delta\rightarrow0.
 $$
  %Without loss of generality, we assume $\check{t}\leq \check{s}$.
  By the properties of $\Psi_{\delta}$, we get that
\begin{eqnarray*}
                 &&\Psi_{\delta}(\check{\gamma}_{\check{t}},\check{\eta}_{\check{s}})%= W'_1(\check{\gamma}_{\check{t}})- {W}'_2(\check{\eta}_{\check{s}})-2\beta(\Upsilon^2(\check{\gamma}_{\check{t}},\hat{{\xi}}_{{\hat{t}}})
%                +\Upsilon^2(\check{\eta}_{\check{s}},\hat{{\xi}}_{{\hat{t}}}))
%                -\frac{1}{\delta}|\check{t}-\check{s}|^2 \\
                  \geq\Psi_{\delta}(\check{\gamma}^{0}_{\check{t}_{0}},\check{\eta}^{0}_{\check{s}_{0}})\geq \sup_{(t,\gamma_t),(s,\eta_s)\in [\hat{t},T]\times \Lambda^{\hat{t}}}\Psi_{\delta}(\gamma_t,\eta_s)-{\delta}\\
                  &\geq&\Psi_{\delta}(\hat{\gamma}_{\hat{t}},\hat{\eta}_{\hat{t}})-{\delta}
                  =\Psi^1(\hat{\gamma}_{\hat{t}},\hat{\eta}_{\hat{t}})-{\delta}%= W_1(\hat{\gamma}_{\hat{t}})- {W}_2(\hat{\eta}_{\hat{t}})-\beta\Upsilon^2(\hat{\gamma}_{\hat{t}},\hat{{\eta}}_{{\hat{t}}})-\delta.
%\end{eqnarray*}
%On the other hand,
%\begin{eqnarray*}
                 %&&\Psi^1(\hat{\gamma}_{\hat{t}},\hat{\eta}_{\hat{t}})
                 \geq \Psi^1(\check{\gamma}_{\check{t},\check{s},A},\check{\eta}_{\check{s}})-{\delta}\\%\geq \Psi_{\delta}(\check{\gamma}_{\check{t},\check{s},A},\check{\eta}_{\check{s}})
                 &=& W'_1(\check{\gamma}_{\check{t},\check{s},A})- {W}'_2(\check{\eta}_{\check{s}})-\beta\Upsilon^2(\check{\gamma}_{\check{t},\check{s},A}-\check{\eta}_{\check{s}})
                 +\varepsilon[\overline{\Upsilon}^2(\check{\gamma}_{\check{t},\check{s},A},\hat{\gamma}_{\hat{t}})
                                     +\overline{\Upsilon}^2(\check{\eta}_{\check{s}},\hat{\eta}_{\hat{t}})]-{\delta}\\
                &=& \Psi_{\delta}(\check{\gamma}_{\check{t}},\check{\eta}_{\check{s}})+\frac{1}{\delta}|\check{t}-\check{s}|^2+W'_1(\check{\gamma}_{\check{t},\check{s},A})-W'_1(\check{\gamma}_{\check{t}})
                +2\beta(\Upsilon^2(\check{\gamma}_{\check{t}},\hat{{\xi}}_{{\hat{t}}})+\Upsilon^2(\check{\eta}_{\check{s}},\hat{{\xi}}_{{\hat{t}}}))-\beta\Upsilon^2(\check{\gamma}_{\check{t},\check{s},A}-\check{\eta}_{\check{s}})\\
                &&+\varepsilon[\overline{\Upsilon}^2(\check{\gamma}_{\check{t},\check{s},A},\hat{\gamma}_{\hat{t}})
                                     +\overline{\Upsilon}^2(\check{\eta}_{\check{s}},\hat{\eta}_{\hat{t}})]-{\delta}.
                %+\frac{\check{s}-\check{t}}{\nu T}(S(\check{\gamma}_{\check{t}})+|\check{\gamma}_{\check{t}}(\check{t})|^2)-\varepsilon|\check{s}-\hat{t}|^2+\varepsilon|\check{t}-\hat{t}|^2\\
               %  &=& \Psi_{\delta}(\check{\gamma}_{\check{t}},\check{\eta}_{\check{s}})+\frac{1}{\delta}|\check{t}-\check{s}|^2+W_1(\check{\gamma}_{\check{t},\check{s},A})-W_1(\check{\gamma}_{\check{t}})
%                +\frac{\check{s}-\check{t}}{\nu T}(S(\check{\gamma}_{\check{t}})+|\check{\gamma}_{\check{t}}(\check{t})|^2)-\varepsilon|\check{s}-\hat{t}|^2+\varepsilon|\check{t}-\hat{t}|^2\\
%                &&+\sum_{i=1}^{\infty}\frac{1}{2^i}[|\check{t}-\hat{t}|^2-|\check{s}-\hat{t}|^2]\\
%                &\geq&\Psi_{\delta}(\check{\gamma}_{\check{t}},\check{\eta}_{\check{s}})+\frac{1}{\delta}|\check{t}-\check{s}|^2+W_1(\check{\gamma}_{\check{t},\check{s},A})-W_1(\check{\gamma}_{\check{t}})
%               -2\varepsilon|\check{s}-\hat{t}||\check{t}-\check{s}|
%                -2\sum_{i=1}^{\infty}\frac{1}{2^i}|\check{s}-\hat{t}||\check{t}-\check{s}|.
\end{eqnarray*}
Noting that, since    $\{e^{tA}, t\geq0\}$ is a $C_0$ contraction semigroup,  for all $\gamma_t\in \Lambda$ and $ s\in [t,T]$,
$$
\Upsilon^2(\gamma_t)=||\gamma_t||_0^2+\frac{|\gamma_t(t)|^4}{||\gamma_t||_0^2}\geq ||\gamma_{t,s,A}||_0^2+\frac{|e^{(s-t)A}\gamma_t(t)|^4}{||\gamma_{t,s,A}||_0^2}=\Upsilon^2(\gamma_{t,s,A})
\geq ||\gamma_{t,s,A}||_0^2=||\gamma_{t}||_0^2\geq\frac{1}{3}\Upsilon^2(\gamma_t).
%\mbox{for all}\   \gamma_t\in \Lambda\ \mbox{and}\  s\in [t,T].
$$
%Then we have
%\begin{eqnarray*}
%                 \Psi^1(\hat{\gamma}_{\hat{t}},\hat{\eta}_{\hat{t}})\geq \Psi_{\delta}(\check{\gamma}_{\check{t}},\check{\eta}_{\check{s}})+\frac{1}{\delta}|\check{t}-\check{s}|^2+W_1(\check{\gamma}_{\check{t},\check{s},A})-W_1(\check{\gamma}_{\check{t}})
%                 +\varepsilon[\overline{\Upsilon}^2(\check{\gamma}_{\check{t},\check{s},A},\hat{\gamma}_{\hat{t}})
%                                     +\overline{\Upsilon}^2(\check{\eta}_{\check{s}},\hat{\eta}_{\hat{t}})].
%                %+\frac{\check{s}-\check{t}}{\nu T}(S(\check{\gamma}_{\check{t}})+|\check{\gamma}_{\check{t}}(\check{t})|^2)-\varepsilon|\check{s}-\hat{t}|^2+\varepsilon|\check{t}-\hat{t}|^2\\
%               %  &=& \Psi_{\delta}(\check{\gamma}_{\check{t}},\check{\eta}_{\check{s}})+\frac{1}{\delta}|\check{t}-\check{s}|^2+W_1(\check{\gamma}_{\check{t},\check{s},A})-W_1(\check{\gamma}_{\check{t}})
%%                +\frac{\check{s}-\check{t}}{\nu T}(S(\check{\gamma}_{\check{t}})+|\check{\gamma}_{\check{t}}(\check{t})|^2)-\varepsilon|\check{s}-\hat{t}|^2+\varepsilon|\check{t}-\hat{t}|^2\\
%%                &&+\sum_{i=1}^{\infty}\frac{1}{2^i}[|\check{t}-\hat{t}|^2-|\check{s}-\hat{t}|^2]\\
%%                &\geq&\Psi_{\delta}(\check{\gamma}_{\check{t}},\check{\eta}_{\check{s}})+\frac{1}{\delta}|\check{t}-\check{s}|^2+W_1(\check{\gamma}_{\check{t},\check{s},A})-W_1(\check{\gamma}_{\check{t}})
%%               -2\varepsilon|\check{s}-\hat{t}||\check{t}-\check{s}|
%%                -2\sum_{i=1}^{\infty}\frac{1}{2^i}|\check{s}-\hat{t}||\check{t}-\check{s}|.
%\end{eqnarray*}
Letting $\delta\rightarrow0$, we obtain that
$$
\frac{1}{\delta}|\check{t}-\check{s}|^2+\varepsilon[\overline{\Upsilon}^2(\check{\gamma}_{\check{t},\check{s},A},\hat{\gamma}_{\hat{t}})
                                     +\overline{\Upsilon}^2(\check{\eta}_{\check{s}},\hat{\eta}_{\hat{t}})]\rightarrow0\ \ \mbox{as}\  \ \delta\rightarrow0,
$$
Notice that
$$
{\Upsilon}^2(\check{\gamma}_{\check{t},\check{s},A},\hat{\gamma}_{\hat{t}})\geq ||\check{\gamma}_{\check{t},\check{s},A}-\hat{\gamma}_{\hat{t},\check{s},A}||_0^2
=||\check{\gamma}_{\check{t}}-\hat{\gamma}_{\hat{t},\check{t},A}||_0^2\geq \frac{1}{3}{\Upsilon}^2(\check{\gamma}_{\check{t}},\hat{\gamma}_{\hat{t}}),
$$
%and
%$$
% W'_1(\check{\gamma}_{\check{t},\check{s}})- {W}'_2(\check{\eta}_{\check{s}})-2\beta(\Upsilon^2(\check{\gamma}_{\check{t}},\hat{{\xi}}_{{\hat{t}}})
%                +\Upsilon^2(\check{\eta}_{\check{s}},\hat{{\xi}}_{{\hat{t}}}))\rightarrow\Psi(\hat{\gamma}_{\hat{t}},\hat{\eta}_{\hat{t}})\ \ \mbox{as}\  \ \delta\rightarrow0.
%$$
  we get that (\ref{4.23}) holds true.  The proof is now complete. \ \ $\Box$
\section{Appendix}  \label{RDS}
%%%%%%%%%%%%%%%%%%%%%%%%%%%%%%%%%%%%%%%%%%%%%%%%%%%%%%%%%%%%%%%%%

\par
In this Appendix, we prove $(\hat{\Lambda}^t, d_{\infty})$ and $({\Lambda}^t, d_{\infty})$ are two complete metric spaces.
\begin{lemma}\label{lemma2.1111}
 $({\Lambda}^t, d_{\infty})$ and $(\hat{\Lambda}^t, d_{\infty})$ are two complete metric spaces for every $t\in [0,T)$. %(see Section 2.1.1 in \cite{tang1}).
\end{lemma}
\par
{\bf  Proof}. \ \
Assume $\{\gamma^n_{t_n}\}_{n\geq0}$ is a cauchy sequence in $(\hat{\Lambda}^t, d_{\infty})$, then for any $\varepsilon>0$,  there exists $N(\varepsilon)>0$ such that, for all $m,n\geq N(\varepsilon)$, we have
\begin{eqnarray*}\label{jialemma111000}
d_\infty(\gamma^n_{t_n},\gamma^m_{t_m})=|t_n-t_m|+\sup_{0\leq s\leq T}|e^{((s-t_n)\vee0)A}\gamma^n_{t_n}(s\wedge t_n)-e^{((s-t_m)\vee0)A}\gamma^m_{t_m}(s\wedge t_m)|<\varepsilon.
\end{eqnarray*}
Therefore, there exists  $\hat{t}\in [t,T]$ such that
%\begin{eqnarray}\label{jialemma2}
$\lim_{n\rightarrow\infty}t_n=\hat{t}.$
%\end{eqnarray}
 Moreover, for all $s\in [0,T]$,
\begin{eqnarray}\label{jialemma}
    |e^{((s-t_n)\vee0)A}\gamma^n_{t_n}(s\wedge t_n)-e^{((s-t_m)\vee0)A}\gamma^m_{t_m}(s\wedge t_m)|<\varepsilon, \  (\forall m,n\geq N(\varepsilon)).
\end{eqnarray}
For fixed $s\in [0,T]$, we see that $\{e^{((s-t_n)\vee0)A}\gamma^n_{t_n}(t_n\wedge s)\}$ is a cauchy sequence, thereby the limit $\lim_{n\rightarrow\infty}e^{((s-t_n)\vee0)A}\gamma^n_{t_n}(t_n\wedge s)$  exists and denoted by $\gamma_{T}(s)$. Letting $m\rightarrow\infty$ in (\ref{jialemma}), we obtain that
\begin{eqnarray*}
    |\gamma_{T}(s)-e^{((s-t_n)\vee0)A}\gamma^n_{t_n}(s\wedge t_n)|\leq\varepsilon, \  (\forall  \ n\geq N(\varepsilon)).
\end{eqnarray*}
 Taking the  supremum over $s\in
             [0,T]$, we get
             \begin{eqnarray}\label{jialemma3}
    \sup_{s\in[0,T]}|\gamma_{T}(s)-e^{((s-t_n)\vee0)A}\gamma^n_{t_n}(s\wedge t_n)|\leq\varepsilon, \  (\forall \ n\geq N(\varepsilon)).
\end{eqnarray}
We claim that $\gamma_{T}(s)=e^{(s-\hat{t})A}\gamma_T(\hat{t})$  for all $s\in(\hat{t},T]$. In fact, if there exists a subsequence $\{t_{n_l}\}_{l\geq0}$ of $\{t_{n}\}_{n\geq0}$ such that $\{t_{n_l}\}_{l\geq0}\leq\hat{t}$, then we have, for every
 $s\in(\hat{t},T]$,
\begin{eqnarray*}
                   \gamma_T(s)&=&\lim_{n\rightarrow\infty}e^{((s-t_n)\vee0)A}\gamma^n_{t_n}(s\wedge t_n)=\lim_{l\rightarrow\infty}e^{(s-t_{n_l})A}\gamma^{n_l}_{t_{n_l}}(t_{n_l})=\lim_{l\rightarrow\infty}e^{(s-\hat{t})A}e^{(\hat{t}-t_{n_l})A}\gamma^{n_l}_{t_{n_l}}(t_{n_l}\wedge \hat{t})\\
                   &=&e^{((s-\hat{t})\vee 0)A}\lim_{n\rightarrow\infty}e^{(\hat{t}-t_{n})A}\gamma^{n}_{t_{n}}(t_{n}\wedge \hat{t})=e^{(s-\hat{t})A}\gamma_T(\hat{t}).
\end{eqnarray*}
 Otherwise, we may assume $\{t_{n}\}_{n\geq0}>\hat{t}$.  For all $s\in (\hat{t},T]$, we can let $m>n$ be large enough such that $t_m\leq s\wedge t_n$, and  letting  $s=t_m$ in (\ref{jialemma}),
 $$
  |e^{(s-t_m)A}\gamma^n_{t_n}(t_m)-e^{(s-t_m)A}\gamma^m_{t_m}(t_m)|\leq M_1|\gamma^n_{t_n}(t_m)-\gamma^m_{t_m}(t_m)|<M_1\varepsilon,
 $$
 Letting  $m\rightarrow\infty$, we obtain
 $$
  |e^{(s-\hat{t})A}\gamma^n_{t_n}(\hat{t})-\gamma_T(s)|\leq M_1\varepsilon.
 $$
 Letting $n\rightarrow\infty$, we have
 $$
 |e^{(s-\hat{t})A}\gamma_{T}(\hat{t})-\gamma_T(s)|\leq M_1\varepsilon,\ \mbox{for all}\ \ s\in(\hat{t},T].
 $$
 %
%
% Then, by the definition of $\gamma_T$, we have that, for $s>\hat{t}$,
% $$
% \lim_{n\rightarrow\infty}[|(e^{(s-\hat{t})A}-e^{(s-t_n)A})\gamma_T(\hat{t})|+|e^{(s-t_n)A}(\gamma_T(\hat{t})-\gamma^n_{t_n}(t_n))|+|e^{(s-t_n)A}\gamma^n_{t_n}(t_n)-\gamma_T(s)|]=0.
% $$
% Thus,
% \begin{eqnarray*}
%                        |e^{(s-\hat{t})A}\gamma_T(\hat{t})-\gamma_T(s)|=0.
% \end{eqnarray*}
% $$
% \lim_{n\rightarrow\infty}[|e^{(s-\hat{t})A}(\gamma_T(\hat{t})-\gamma^n_{t_n}(t_n))|+|e^{(s-t_n)A}\gamma^n_{t_n}(t_n)-\gamma_T(s)|]=0;
% $$
%\begin{eqnarray*}
% &&\lim_{n\rightarrow\infty}|(e^{(s-\hat{t})A}-e^{(s-t_n)A})\gamma^n_{t_n}(t_n)|\\
% &\leq& \lim_{n\rightarrow\infty}[|e^{(s-\hat{t})A}(\gamma^n_{t_n}(t_n)-\gamma_T(\hat{t}))|+|(e^{(s-\hat{t})A}-e^{(s-t_n)A})\gamma_T(\hat{t})|+|e^{(s-t_n)A}(\gamma_T(\hat{t})-\gamma^n_{t_n}(t_n))|]=0.
%\end{eqnarray*}
%Thus,
% \begin{eqnarray*}
%                        &&|e^{(s-\hat{t})A}\gamma_T(\hat{t})-\gamma_T(s)|\\
%                        &\leq& \lim_{n\rightarrow\infty}[|e^{(s-\hat{t})A}(\gamma_T(\hat{t})-\gamma^n_{t_n}(t_n))|+|(e^{(s-\hat{t})A}-e^{(s-t_n)A})\gamma^n_{t_n}(t_n)|+|e^{(s-t_n)A}\gamma^n_{t_n}(t_n)-\gamma_T(s)|]=0.
% \end{eqnarray*}
Then, by (\ref{jialemma3}) we obtain
\begin{eqnarray}\label{jialemma23}
d_\infty(\eta_{\hat{t}},\gamma^n_{t_n})\rightarrow0 \ \mbox{as} \  n\rightarrow\infty.
\end{eqnarray}
Here we let $ \eta_{\hat{t}}$ denote $\gamma_T|_{[0,\hat{t}]}$.
Now we prove $\eta_{\hat{t}}\in \hat{\Lambda}^t$. First, we prove $\eta_{\hat{t}}$ is right-continuous.
For every $0\leq s<\hat{t}$ and $0<\delta< \hat{t}-s$, we have
$$
|\eta_{\hat{t}}(s+\delta)-\eta_{\hat{t}}(s)|\leq |\gamma_{T}(s+\delta)-\gamma^n_{t_n}((s+\delta)\wedge t_n)|+|\gamma^n_{t_n}((s+\delta)\wedge t_n)-\gamma^n_{t_n}(s\wedge t_n)|+|\gamma^n_{t_n}(s\wedge t_n)-\gamma_{T}(s)|.
$$
For every $\varepsilon>0$, by (\ref{jialemma3}), there exists $n>0$ be large enough such that
$$
|\gamma_{T}(s+\delta)-\gamma^n_{t_n}((s+\delta)\wedge t_n)|+|\gamma^n_{t_n}(s\wedge t_n)-\gamma_{T}(s)|< \frac{\varepsilon}{2}.
$$
For the fixed $n$, since $\gamma^n_{t_n}\in \hat{\Lambda}^t$, there exists a constant $0<\Delta< \hat{t}-s$ such that, for all $0\leq\delta<\Delta$,
$$
  |\gamma^n_{t_n}((s+\delta)\wedge t_n)-\gamma^n_{t_n}(s\wedge t_n)|< \frac{\varepsilon}{2}.
$$
Then $
|\eta_{\hat{t}}(s+\delta)-\eta_{\hat{t}}(s)|< \varepsilon
$ for all $0\leq\delta<\Delta$. Next, let us prove  $\eta_{\hat{t}}$  has left limit in $(0,\hat{t}]$. For every $0< s\leq\hat{t}$ and $0\leq s_1,s_2<s$, we have
$$
|\eta_{\hat{t}}(s_1)-\eta_{\hat{t}}(s_2)|\leq |\gamma_{T}(s_1)-\gamma^n_{t_n}(s_1\wedge t_n)|+|\gamma_{T}(s_2)-\gamma^n_{t_n}(s_2\wedge t_n)|+|\gamma^n_{t_n}(s_1\wedge t_n)-\gamma^n_{t_n}(s_2\wedge t_n)|.
$$
For every $\varepsilon>0$, by (\ref{jialemma3}), there exists $n>0$, which is independent of $s_1,s_2$, be large enough such that
$$
|\gamma_{T}(s_1)-\gamma^n_{t_n}(s_1\wedge t_n)|+|\gamma_{T}(s_2)-\gamma^n_{t_n}(s_2\wedge t_n)|< \frac{\varepsilon}{2}.
$$
For the fixed $n$,% if $t_n<s$, we can let $\Delta>0$ be small enough such that $t_n<s-\Delta$, then for all  $s_1,s_2\in [s-\Delta,s)$
%$$
%|\gamma^n_{t_n}(s_1\wedge t_n)-\gamma^n_{t_n}(s_2\wedge t_n)|=|\gamma^n_{t_n}(t_n)-\gamma^n_{t_n}(t_n)|=0;
%$$
%if $t_n\geq s$,
  since $\gamma^n_{t_n}\in \hat{\Lambda}^t$, there exists a constant $\Delta>0$ such that, for all $s_1,s_2\in [s-\Delta,s)$,
$$
|\gamma^n_{t_n}(s_1\wedge t_n)-\gamma^n_{t_n}(s_2\wedge t_n)|
%=|\gamma^n_{t_n}(s_1)-\gamma^n_{t_n}(s_2)|
<\frac{\varepsilon}{2}.
$$
Then there exists a constant $\Delta>0$ such that $
|\eta_{\hat{t}}(s_1)-\eta_{\hat{t}}(s_2)|< \varepsilon
$ for all  $s_1,s_2\in [s-\Delta,s)$.
\par
         Finally, we prove  $\eta_{\hat{t}}\in {\Lambda}^t$ if $\{\gamma^n_{t_n}\}_{n\geq0}$ is a cauchy sequence in $({\Lambda}^t, d_{\infty})$.
        By the similar proof process of right continuous above, we get that $\eta_{\hat{t}}$ is continuous in $[0,\hat{t})$.
        %
%         For every $0\leq t, s<\hat{t}$, we have
%$$
%|\eta_{\hat{t}}(t)-\eta_{\hat{t}}(s)|\leq |\gamma_{T}(t)-\gamma^n_{t_n}(t\wedge t_n)|+|\gamma^n_{t_n}(t\wedge t_n)-\gamma^n_{t_n}(s\wedge t_n)|+|\gamma^n_{t_n}(s\wedge t_n)-\gamma_{T}(s)|.
%$$
%For every $\varepsilon>0$, by (\ref{jialemma3}), there exists $n>0$ be large enough such that
%$$
%|\gamma_{T}(t)-\gamma^n_{t_n}(t\wedge t_n)|+|\gamma^n_{t_n}(s\wedge t_n)-\gamma_{T}(s)|< \frac{\varepsilon}{2}.
%$$
%For the fixed $n$, since $\gamma^n_{t_n}\in {\Lambda}^t$, there exists a constant $\Delta>0$ such that, for all $|t-s|<\Delta$,
%$$
%  |\gamma^n_{t_n}(t\wedge t_n)-\gamma^n_{t_n}(s\wedge t_n)|< \frac{\varepsilon}{2}.
%$$
%Then $
%|\eta_{\hat{t}}(t)-\eta_{\hat{t}}(s)|< \varepsilon
%$ for all $|t-s|<\Delta$.
 Now we have to prove that $\eta_{\hat{t}}$ is left continuous at $\hat{t}$. For every $0\leq t<\hat{t}$, we have, if $t_n\geq \hat{t}$,
$$
|\eta_{\hat{t}}(t)-\eta_{\hat{t}}(\hat{t})|\leq |\gamma_{T}(t)-\gamma^n_{t_n}(t\wedge t_n)|+|\gamma^n_{t_n}(t\wedge t_n)-\gamma^n_{t_n}(\hat{t}\wedge t_n)|+|\gamma^n_{t_n}(\hat{t}\wedge t_n)-\gamma_{T}(\hat{t})|.
$$
For every $\varepsilon>0$, by  the similar process above, there exists a constant $\Delta>0$ such that, for all $|t-\hat{t}|<\Delta$,
$$
  |\eta_{\hat{t}}(t)-\eta_{\hat{t}}(\hat{t})|< \varepsilon.
$$
If $t_n<\hat{t}$,
\begin{eqnarray*}
|\eta_{\hat{t}}(t)-\eta_{\hat{t}}(\hat{t})|&\leq& |\gamma_{T}(t)-e^{(({t}-t_n)\vee0)A}\gamma^n_{t_n}(t\wedge t_n)|+|e^{(({t}-t_n)\vee0)A}\gamma^n_{t_n}(t\wedge t_n)-e^{(\hat{t}-t_n)A}\gamma^n_{t_n}( t_n)|\\
&&+|e^{(\hat{t}-t_n)A}\gamma^n_{t_n}(t_n)-\gamma_{T}(\hat{t})|.
\end{eqnarray*}
For every $\varepsilon>0$, by (\ref{jialemma3}), there exists $n>0$ be large enough such that
$$
|\gamma_{T}(t)-e^{(({t}-t_n)\vee0)A}\gamma^n_{t_n}(t\wedge t_n)|+|e^{(\hat{t}-t_n)A}\gamma^n_{t_n}(t_n)-\gamma_{T}(\hat{t})|< \frac{\varepsilon}{2}.
$$
For the fixed $n$, since  $\{e^{tA}, t\geq0\}$ is a $C_0$
                 semigroup, there exists a constant $0<\Delta\leq |\hat{t}-t_n|$ such that, for all $|t-\hat{t}|<\Delta\leq |\hat{t}-t_n|$,
$$
  |e^{(({t}-t_n)\vee0)A}\gamma^n_{t_n}(t\wedge t_n)-e^{(\hat{t}-t_n)A}\gamma^n_{t_n}( t_n)|=|(e^{(\hat{t}-t)A})-I)e^{({t}-t_n)A}\gamma^n_{t_n}(t_n)|< \frac{\varepsilon}{2}.
$$
Then $
|\eta_{\hat{t}}(t)-\eta_{\hat{t}}(\hat{t})|< \varepsilon
$ for all $|t-\hat{t}|<\Delta$.
The proof is now complete. \ \ $\Box$

\par

\end{document}